\documentclass[a4paper,11pt, reqno]{amsart}
\usepackage[margin=2.5cm]{geometry}
\usepackage{hyperref}
\usepackage{amsmath,amsfonts}
\usepackage{amsthm}
\usepackage[authoryear]{natbib}
\usepackage{enumitem}% http://ctan.org/pkg/enumitem

\usepackage{amssymb}

\let\emptyset\varnothing
\usepackage{multirow}
\usepackage{float}
\usepackage{subcaption}
\usepackage{color}
\usepackage{adjustbox}
\usepackage{booktabs}
\usepackage{soul}
\def\R{\mathbb{R}}

\def\B{\mathbb{B}}

\def\A{\mathcal{A}}
\def\B{\mathbb{B}}

\def\Cov{\mathcal{C}}
\newcommand{\dsum}{\displaystyle\sum}

\newtheorem{theorem}{Theorem}

\newtheorem{lem}[theorem]{Lemma}
\newtheorem{prop}[theorem]{Proposition}
\newtheorem{example}{Example}
\newtheorem{remark}{Remark}

\allowdisplaybreaks

\let\origmaketitle\maketitle
\def\maketitle{
  \begingroup
  \def\uppercasenonmath##1{} % this disables uppercasing title
  \let\MakeUppercase\relax % this disables uppercasing authors
  \origmaketitle
  \endgroup
}

\begin{document}
\title[]{\Large Multitype Maximal Covering Location Problems:\\Hybridizing discrete and continuous problems}

\author[V. Blanco, R. G\'aquez \MakeLowercase{and} F. Salanha-da-Gama]{{\large V\'ictor Blanco$^\dagger$, Ricardo G\'azquez$^\dagger$ and  Francisco Saldanha-da-Gama$^\ddagger$}\medskip\\
$^\dagger$Institute of Mathematics (IMAG), Universidad de Granada\\
$^\ddagger$Centro de Matemática, Aplicações Fundamentais e Investigação Operacional, Faculdade de Ciências da Universidade de Lisboa\\
\texttt{vblanco@ugr.es}, \texttt{rgazquez@us.es}, \texttt{faconceicao@fc.ul.pt}
}

\date{\today}

\maketitle 
\begin{abstract}
	This paper introduces a general modeling framework for a multi-type maximal covering location problem in which the position of facilities in different metric spaces are simultaneously decided to maximize the demand generated by a set of points. From the need of intertwining location decisions in discrete and in continuous sets, a general hybridized problem is considered in which some types of facilities are to be located in finite sets and the others in continuous metric spaces.
	A natural non-linear model is proposed for which an integer linear programming reformulation is derived.
	A branch-and-cut algorithm is developed for better tackling the problem.
	The study proceeds considering the particular case in which the continuous facilities are to be located in the Euclidean plane. 
	In this case, taking advantage from some geometrical properties it is possible to propose an alternative integer linear programming model. 
	The results of an extensive battery of computational experiments performed to assess the methodological contribution of this work is reported on.
	The data consists of up to 920 demand nodes using real geographical and demographic data. 
	\end{abstract}


\section{Introduction}
Among the most popular problems in Location Analysis are those in which a user can receive a service in case it is located close enough to an open facility providing it. 
These problems are usually called \textit{Covering Location problems}. 
Two different paradigms have been considered when classifying this type of problems. 
The first considers  a cost-oriented objective and the main goal is to satisfy the demand of \textit{all} the users by minimizing the setup costs of the facilities. 
These problems are referred to in the literature as \textit{Set Covering Problems}, and were mathematically introduced by \citet{toregas1971location}.
The most popular problem in this family is the $p$-center problem~\citep{Hakimi:1965}.
On the other hand, the second paradigm assumes the existence of a budget for opening facilities and the goal is to accommodate it to satisfy as much demand of the users as possible. 
These problems belong to the family of \textit{Maximal Covering Location Problems} (MCLP) that have attracted the attention of many researchers since its introduction by \citet{ChurchReVelle:1974}, both because its practical interest in different disciplines see \citep[see][]{Chung:1986} and the mathematical challenges it  poses.
The interested readers are referred to the surveys by \citet{Church&Murray:2018} and \citet{GarciaMarin:2019} for further details on covering location problems. 

As in most facility location problems, the MCLP have been analyzed within three main frameworks attending to the space where the facilities can be located: discrete~\citep[see][]{arana2020fuzzy,cordeau2019benders}, network  \citep[see][]{baldomero2021minmax,berman2016covering,church1979location}  and continuous~\citep[see][]{blanco2021continuous,drezner2019planar,he2016mean,tedeschi2021new}.
The applicability of each of these frameworks mainly depends on the nature of the facilities to be located and the decision maker's knowledge in terms of the available information.
While in a discrete setting the centers to open are selected from a finite set of possibilities, in continuous location the facilities are allowed to be located in the whole underlying space. 
In network problems the centers can be located anywhere in a network.
Like for many minsum facility location problems on networks, in the case of covering problems often a finite dominating set can be identified which allows casting the problem in a discrete setting giving great prominence to the research done in the latter. 

Discrete facility location is mostly applied to locate physical services (e.g., ATMs, stores, hospitals); continuous location is known to be more adequate when locating telecommunication services (e.g. routers, surveillance equipment, alarms, sensors, etc.) where more flexibility is allowed for their optimal positions.

A feature shared by most of the exiting literature focusing the MCLP concerns the existence of a single type of facility.
However, in practice, this may not be the case.
If not by other reasons, the progressive technology development often calls for older equipment that is still operational to be used together with more recent one.
Another possibility emerges when two technologies can be looked at as complementing each other.
For instance, when locating equipment for early fire detection in forests, surveillance facilities requiring human resources operating them may be complemented with equipment such as remotely controlled cameras to ensure a better coverage of the area of interest.
When facilities can be installed in different phases (e.g. multi-period facility location) the facilities to be located in each phase can also be looked at as belonging to a different group (that we still call type) of facilities.

In this paper we investigate maximal covering location problems with multiple facility types. 
We assume that the number of facilities of each type to be located is known beforehand and that each type of facilities is characterized by the shape of their coverage areas and the metric space from which they are selected.
A plan is to be devised for a multi-stage process with each stage corresponding to installing one type of facility.
The facilities opened in each stage perform the same tasks and thus complement the facilities installed in previous stages.
In turn, they will be complemented by the facilities to be located in the subsequent stages.
We show that instead of making sequential decisions for each facility type, coverage gains can be achieved by making an integrated decision involving all facility types.

We start by presenting a general formulation for the problem.
Afterwards, motivated by some practical settings we investigate the hybridization of discrete and continuous facility location.
We consider that several types of facilities are to be selected in finite sets of possibilities (one for each type) whereas the other types of facilities can be located continuously in the whole space.
%We specify the above hybridization for two-dimensional Euclidean spaces.

By considering a hybrid setting it becomes possible to take advantage from choosing some services in finite sets of pre-specified preferred locations and then deciding flexible positions of the servers in the whole space.
This setting can be useful, for instance, in telecommunications networks with a certain number of the servers (sensors, antennas, routers, etc.) being located inside adequately prepared infrastructures (buildings, offices, air-conditioned cabins, roofs, etc.) and additional servers being located at any place in the given space.
The goal is of course to capture/cover as much demand as possible no matter the equipment doing it. 
The continuous facilities can be looked at as a set of servers to be located in the future and that must complement the equipment located in a discrete setting. 
To decide the positions of the centers, one could proceed by first locating the initial centers (in a fully discrete framework) maximizing the covered demand and then locate the future centers (in a fully continuous framework), maximizing the covered demand of the customers that are still uncovered by the initial servers.  
Although allowing the application of well-known existing tools in the context discrete and continuous maximal covering location, this procedure may easily lead to sub-optimal solutions: a better planning (i.e., covering more demand) can be obtained by considering an integrated approach which is what we propose. 

The literature is quite rich when it comes to considering multi-type facility location problems.
Nevertheless, more often than not, we are led to problems stemming from logistics or telecommunications applications in which a multi-layered or a hierarchical facility structure is to be setup. 
In such a case, each facility type lies within a specific layer of the network or in the hierarchy and has a specialized function.
The reader can refer to \citet{ContrerasOrtiz-Astorquiza:2019} and \citet{HechmannNickel:2019} as well as to the references there in for overviews on many such problems.

In the current work, we are concerned with facilities that provide the same service, and thus can be used to complement each other although having some different characteristics.
\citet{WuZhangZhang:2006} investigated one such problem in the context of capacitated facility location.
In that paper, general setup costs are considered that depend on the size and location of the facility.
The problem lies in the context of fixed-charge facility location \citep{FernancezLandete:2019}.

\citet{Mesa:1991} investigated several multi-period problems on networks.
In particular, the author introduced the so-called absolute multi-period $(\alpha_1,\dots,\alpha_{|T|})-$median problem where $T$ stands for the number of periods in the planning horizon. 
This is possibly the first multi-period extension of the network $p$-median. 
We can look at the facilities to be located in each period as being of different types. 
Unlike we are considering in the current paper, some facilities may just replace others (the former are closed when the latter are opened) and the location space is the same for all facilities.

\citet{BermanDrezner:2008} consider a two-type discrete facility location problem.
This is a $p$-median problem under uncertainty consisting of locating $p$ initial facilities plus an uncertain number of extra additional ones. 
In both cases, the potential facilities to open belong to a finite set that coincides with the demand points and thus the problem is cast as a stochastic discrete $p$-median problem.
In the current paper we assume that we know beforehand the number of each type of facility to locate.
Furthermore, we consider specific location spaces according to the different types of facilities.

\citet{HeynsVanVuuren:2018} investigate a problem in which multiple types of facilities can be located in specific zones identified beforehand.
Type-specific location requirements are assumed for the facilities.
In each zone a finite set of candidate locations for the facilities are assumed, i.e., a pure discrete facility location setting is adopted.
All facility types can in principle be located in all zones.

Considering also a finite set for locating the facilities, we find works considering a hierarchy between the facility-types in line to the models discussed by \citet{ContrerasOrtiz-Astorquiza:2019}, i.e., the facilities in a higher level extend the service provided by the facilities in lower levels.
\citet{moore1982hierarchical} are possibly the first authors introducing such type of problem.
In each potential facility location one must decide the type of facility to locate (if some).
Other works deepening this type of analysis include \citet{espejo2003dual}, \citet{ratick2009application}, and \citet{Xia&Yin&Xu&Xie&Shao:2009}.
More recently, \citet{KucukaydinAras:2020} also investigate a multi-type discrete facility location problem but they consider so-called consumer preference.
Each demand point has a preference for one facility type.
The facilities are to be located in such a way that an optimal coverage in terms of the consumer preference is achieved.

In the current paper we go beyond the existing literature by proposing a general modeling framework for a multi-type maximal covering location problem.
We do not restricted the problem to a discrete setting.
Instead, we assume some location space for each facility type.
The general framework proposed is motivated by some applications calling for hybridizing discrete and continuous facility location problems. 
For this reason, we deepen the study by considering that hybridized case: several types of facilities are to be located in a finite set of possible locations with their service being complemented by other facility types that can be located anywhere in underlying continuous space.
We propose a `natural' non-linear model that nonetheless, raises some computational difficulties. 
For this reason we also develop an integer linear model.
Afterwards we focus on on the Euclidean plane.
This allows using other types of modeling frameworks that we also investigate.
We report on a series results obtained from a series of computational tests performed to assess the different models proposed.
Real geographical data is consider in these tests.

The remainder of the paper is organized as follows. In Section~\ref{sec:generalProblem} the investigated problem is detailed and a general mathematical model is introduced. Section~\ref{sec:hybridizedDiscCont} specializes the general modeling framework to a hybridized discrete-continuous setting. Afterwards, in Section~\ref{sec:EuclideanPlane} we focuses on the the Euclidean plane provided additional insights in this case. The results of the computational tests performed to asses different developments proposed in the paper are reported in Section~\ref{sec:ComputationalExperiments}.
The paper end with an overview of the work done and some suggestions for further research.

\section{The Multitype Maximal Covering Location Problem}\label{sec:generalProblem}

Consider a finite set $\A=\{a_1,\dots, a_n\}$ of demand points in $\R^d$ indexed in set $N=\{1,\dots,n\}$, each of which with a weight given by a non-negative value $c_i$ representing the demand of node $a_i$, for all $i \in N$. 
Throughout the paper we often call a demand point interchangeably by the node $a_i$ or by the index $i$.

Let us assume that there is a finite set of facility types, indexed in a set $\mathcal{T}=\{1,\dots,T\}$.
A facility of type $t \in \mathcal{T}$ can be located in some metric space that we denote by $\mathbb{S}^{(t)} \subseteq \R^d$.
Consider a distance function of interest, say $\|\cdot\|^{(t)}$, in $\mathbb{S}^{(t)}$.
Also for a facility of type $t \in \mathcal{T}$ we consider a coverage radius, say $\rho(t)$, $t \in \mathcal{T}$. 
Given one such facility, we say that node $a_i$, $i \in N$, is \emph{covered by the facility} if the distance between $a_i$ and the facility does not exceed $\rho(t)$.
The node $a_i$ is said to be \textit{covered} if there is at least one open facility (no matter its type) covering it. 
For any finite subset of open facilities of type $t\in \mathcal{T}$, $\mathcal{X}^{(t)} \subseteq \mathbb{S}^{(t)}$, we denote by $\Cov(\mathcal{X}^{(t)}) \subseteq N$ the indices of the nodes in $\mathcal{A}$ covered by at least one point in $\mathcal{X}^{(t)}$, i.e.,
$$
\Cov(\mathcal{X}^{(t)}) = \{i \in N: \|a_i-b\|^{(t)} \leq \rho(t), \text{for some $b \in \mathcal{X}^{(t)}$} \}.
$$
%Thus, given 
Given $\mathbf{p}=(p_1, \ldots, p_T) \in N \times\stackrel{T}{\cdots}\times N$, the problem that we call the $\mathbf{p}$-Multitype Maximal Covering Location Problem ($\mathbf{p}$-MTMCLP, for short) seeks to locate $p_t$ facilities of type $t \in \mathcal{T}$ so that the covered demand is maximized.

With the above notation, the $\mathbf{p}$-MTMCLP can be formally formulated as the following optimization problem:
\begin{equation}\label{p1...pTMCLP}\tag{$\mathbf{p}$-{\rm MTMCLP}}
\mathcal{V}(\mathbf{p}) := \max_{\mathcal{X}^{(t)} \subseteq \mathbb{S}^{(t)}, \: |\mathcal{X}^{(t)}|=p_t} \:\: \sum_{i \in \bigcup_{t \in \mathcal{T}} \Cov(\mathcal{X}^{(t)})} c_i.
\end{equation}

Observe that the difference between the different types of facilities to be located is the metric space where the locations are to be found as well as the coverage radii. 
In case the metric spaces coincide for two different types of facilities, one may also consider that they are of the same type and define different coverage radii, resulting in the same model.

The following example illustrates the problem we are investigating. 
Furthermore, it shows that a sequential decision making process in which we locate one type of facility in each step may lead to a sub-optimal solution.
\begin{example}\label{ex:1}
	We randomly generated a set of 50 demand nodes in $[0,1] \times [0,1]$---our set $\A$.
	Let us assume three types of facilities with $\mathbb{S}^{(1)} = \mathcal{A}$ and $\mathbb{S}^{(2)} = \mathbb{S}^{(3)} =\R^2$.
	Additionally, assume that $\|\cdot\|^{(1)} = \|\cdot\|^{(2)}$ are the Euclidean norm and $\|\cdot\|^{(3)}$ is the $\ell_3$-norm.
	Regarding the coverage radii we take $\rho(1) = 0.2$ and $\rho(2)=\rho(3)=0.1$.  
	The number of facilities to open was fixed to $\mathbf{p}=(2,2,1)$. 
	The weights $c_i$ were all set equal to one.
	
	\begin{figure}[ht]
		\begin{subfigure}[b]{.33\linewidth}
			\centering 
			\fbox{\includegraphics[scale=0.4]{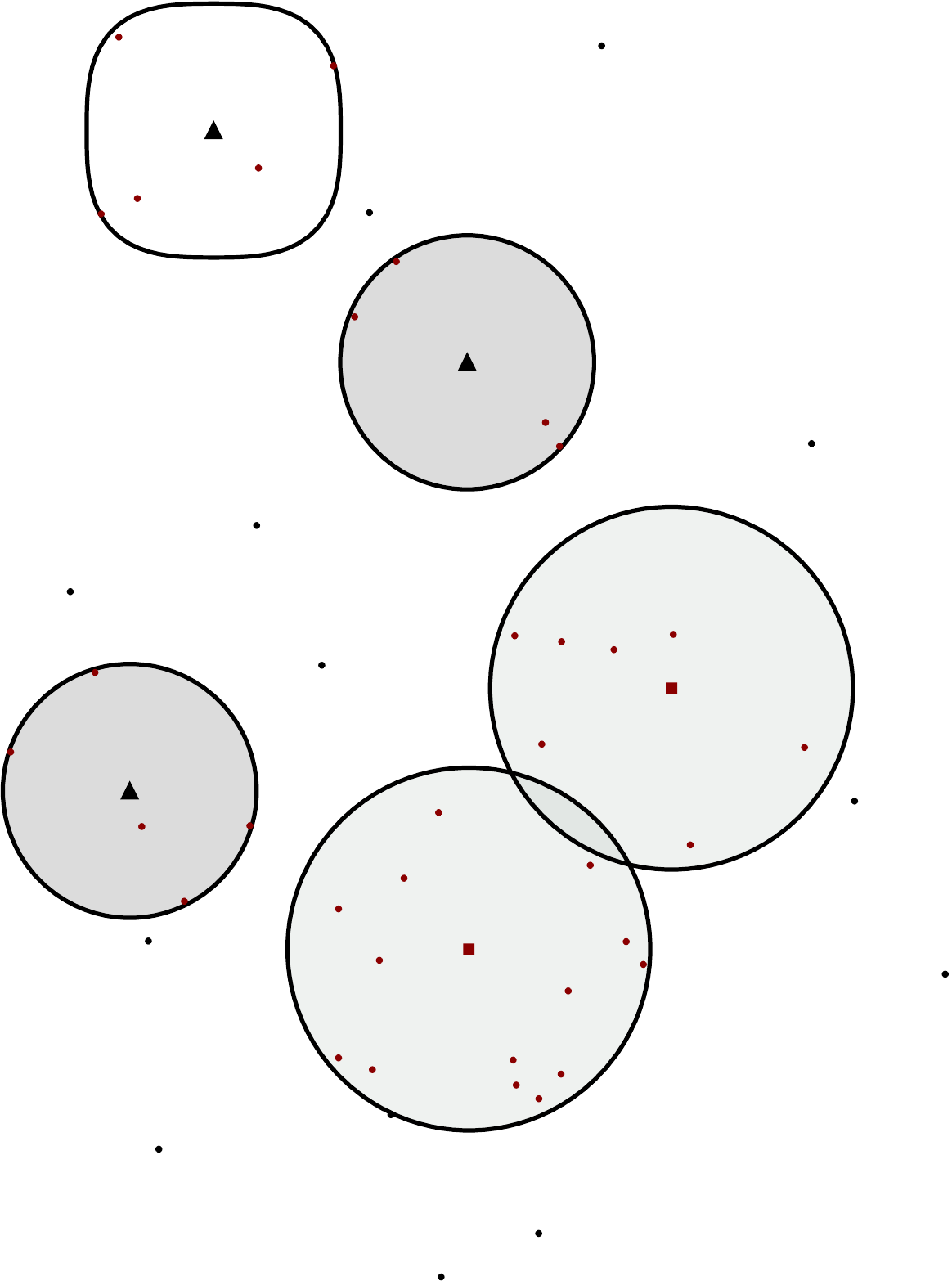}}
			\caption{${(\mathbb{S}^{(1)}\rightarrow\mathbb{S}^{(2)}\&\mathbb{S}^{(3)})}$}\label{fig0a}
		\end{subfigure}~\begin{subfigure}[b]{.33\linewidth}
			\centering 
			\fbox{\includegraphics[scale=0.4]{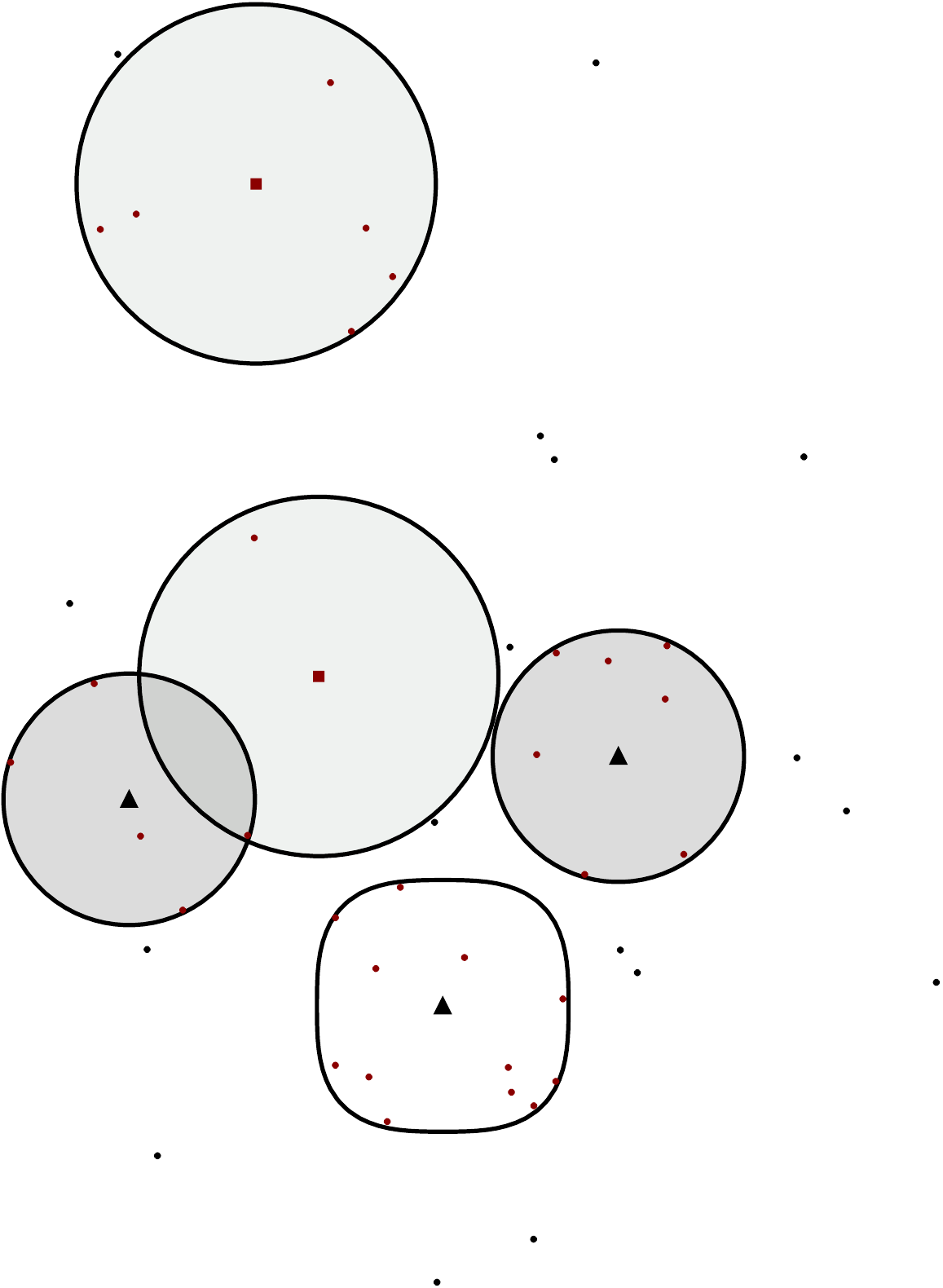}}
			\caption{${(\mathbb{S}^{(2)}\&\mathbb{S}^{(3)}\rightarrow\mathbb{S}^{(1)})}$}\label{fig0b}
		\end{subfigure}~\begin{subfigure}[b]{.33\linewidth}
			\centering 
			\fbox{\includegraphics[scale=0.4]{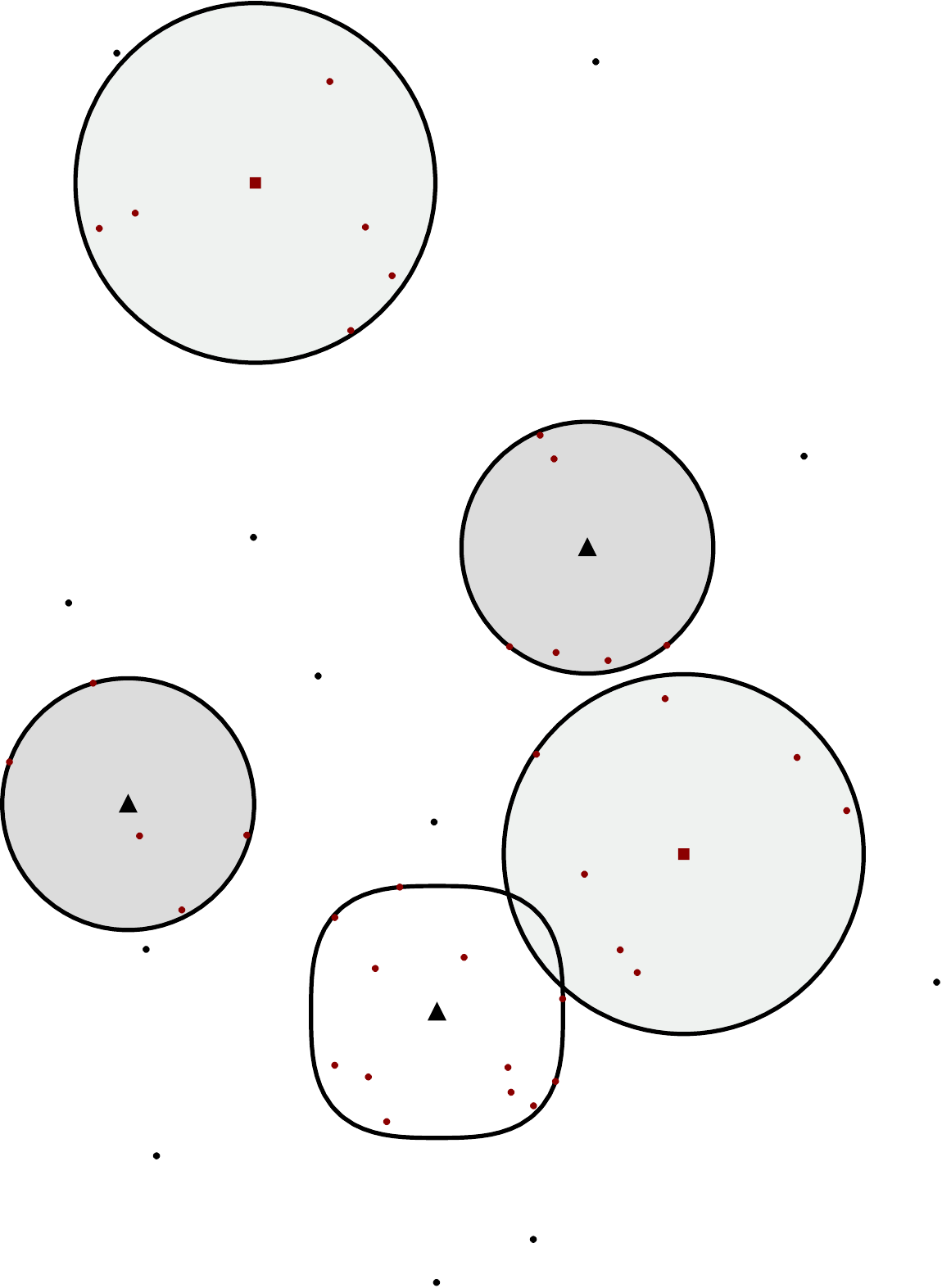}}
			\caption{$\mathbf{p}$-${\rm MTMCLP}$}\label{fig0c}
		\end{subfigure}~
		\caption{Sequential versus integrated decision making.}\label{fig0}
	\end{figure}
	
	First, we look into locating two facilities in $\mathbb{S}^{(1)}$ maximizing the number of covered points in $\A$, and after that we look for the best three additional facilities (two in $\mathbb{S}^{(2)}$ and one in $\mathbb{S}^{(3)}$) maximizing the number of covered points (i.e., covered demand) that were not already covered with the two initial facilities. 
	We denote the resulting solution as ${(\mathbb{S}^{(1)} \rightarrow \mathbb{S}^{(2)}\&\mathbb{S}^{(3)})}$. 
	Using the same instance, a different solution is obtained when we revert the sequence of decisions, i.e., we first look simultaneously for the two + one facilities in $\mathbb{S}^{(2)}$ and $\mathbb{S}^{(3)}$, respectively and after that we seek to find the additional two facilities in $\mathbb{S}^{(1)}$ maximizing the  still not covered demand. The resulting solution is denoted by ${(\mathbb{S}^{(2)}\&\mathbb{S}^{(3)}  \rightarrow\mathbb{S}^{(1)})}$. In a third experiment, we compute the solution using the integrated model \eqref{p1...pTMCLP}.

	The solution ${(\mathbb{S}^{(1)} \rightarrow \mathbb{S}^{(2)}\&\mathbb{S}^{(3)})}$, which is depicted in Figure~\ref{fig0a}, is such that 74\% of the demand nodes are covered. 
	In case of solution ${(\mathbb{S}^{(2)}\&\mathbb{S}^{(3)}  \rightarrow\mathbb{S}^{(1)})}$, depicted in Figure~\ref{fig0b}, this percentage decreases to 66\%. 
	Finally, when we consider the solution obtained using model ($\mathbf{p}$-${\rm MTMCLP}$)---Figure~\ref{fig0c}---we obtain a 76\% demand coverage.  
	These results show that a sequential decision making process (even aggregating some types) may lead to a sub-optimal solution, which gives strength to the integrated modeling framework we are investigating.
	
	It is also worth noticing the change in the shape of the covered area associated with facility to be located is $\mathbb{S}^{(3)}$.
\end{example}

The above model is much general and accommodates many situations as discussed in the introductory section.
In particular the location spaces $\mathbb{S}^{(t)}$, $t \in \mathcal{T}$, may correspond to finite sets, networks, continuous spaces or a combination of all.
\citet{HeynsVanVuuren:2018} and \citet{KucukaydinAras:2020} consider a pure discrete setting.
In the next section we focus on a case that raises some interesting challenges and that corresponds to having some types of facilities to be located in discrete spaces and the others located in continuous ones.

\section{The hybridized discrete-continuous Maximal Covering Location Problem}\label{sec:hybridizedDiscCont}

In this section we 
%analyze 
propose suitable mathematical programming formulations for a wide family of problems in the shape of \eqref{p1...pTMCLP} namely, the one that result when one assumes that the metric spaces $\mathbb{S}^{(t)}$ are either finite sets of points or the entire space $\R^d$.

We keep considering the set of demand points $\mathcal{A}$ already introduced as well as the index set for it elements, $N=\{1,\dots,n\}$. We assume that two main families of types of facilities are to be located, say $T_1$ of type \emph{discrete} and $T_2$ of type \emph{continuous}, so we have a total of $T =  T_1 + T_2$ types of facilities. 
In particular, we consider the set of type indices given by $\mathcal{T}=\mathcal{T}_1\cup \mathcal{T}_2$ where $\mathcal{T}_1 =\{1,\ldots,T_1\}$, and $\mathcal{T}_2=\{T_1+1,\ldots,T\}$.

For the first $T_1$ types of facilities, we consider that $\mathbb{S}^{(t)}=\{b^{(t)}_1,\dots,b^{(t)}_{m^{(t)}}\}$ is a finite set of points in $\R^d$, indexed in set $M^{(t)}=\{1,\dots,m^{(t)}\}$ that have been identified as potential locations for the facilities of type $t$, for $t\in \mathcal{T}_1$. 

As in the general framework, a facility of type $t\in \mathcal{T}$ is endowed with a coverage radius $\rho(t)$. 
However, for the finite location spaces, we can go further in terms of coverage radii specification by assuming facility-dependent radii for the facilities in $\mathbb{S}^{(1)},\dots,\mathbb{S}^{(T_1)}$.
We assume that a facility located at $b^{(t)}_j$ is endowed with a coverage radius equal to $\rho_j(t)$.
In fact, the hybridized discrete-continuous setting we are considering is motivated by some practical applications (e.g. in telecommunication networks planning) an thus it makes sense to consider different coverage areas (typically larger) for the facilities chosen from the pre-specified sets $\mathbb{S}^{(1)},\dots,\mathbb{S}^{(T_1)}$ than those for the extra facilities---facilities located in $\mathbb{S}^{(T_1+1)},\dots,\mathbb{S}^{(T)}$ since the physical infrastructures may be prepared for a better service. 
Moreover, the fact that we know the potential locations in advance allows specifying a coverage radius that is location-specific.
On the other hand, the common radii assumed for the different types of continuous facilities reflect a guaranteed coverage provided by the equipment no matter the point in the space it will end up being located. 

For each $t\in \mathcal{T}$ we denote by $P_t$ the index set for the facilities of type $t$ to be located, i.e., $P_t=\{1, \ldots, p_t\}$.

For $t\in \mathcal{T}_2$, we consider metric spaces $\mathbb{S}^{t} = \R^d$ such that $p_{T_1+1},\ldots, p_T$  locations are to be found respectively, in each of these spaces for installing additional facilities.
These facilities have a coverage radius equal to $\rho(t)$, respectively for type $t\in \mathcal{T}_2$. 
In case one requires the continuous facilities to be located in a specific region of $\R^d$ instead of the entire space, most of the results presented in this paper can be adapted conveniently in case the sets $\mathbb{S}^{(t)}$ are polyhedra or second order cone representable sets, by adding the suitable constraints defining each specific set.

In the following example we illustrate the problem on a three-dimensional instance.
\begin{example}\label{ex:2}
		We randomly generated a set of 50 demand nodes in $[0,1]\times[0,1]\times[0,1]-$our set $\A$. In this case, we consider that $\|\cdot\|^{(1)} = \|\cdot\|^{(2)}$ is the Euclidean norm. Finally, the rest of parameters are selected in the same way: two types of facilities with $\mathbb{S}^{(1)} = \mathcal{A}$ and $\mathbb{S}^{(2)} =\R^3$, the coverage radii we take $\rho(1) = 0.2$ and $\rho(2)=0.1$, the number of facilities to open was fixed to $\mathbf{p}=(2,3)$, and the weights $c_i$ were all set equal to one.
		
		\begin{figure}[H]
				\centering 
				\includegraphics[width=0.6\textwidth]{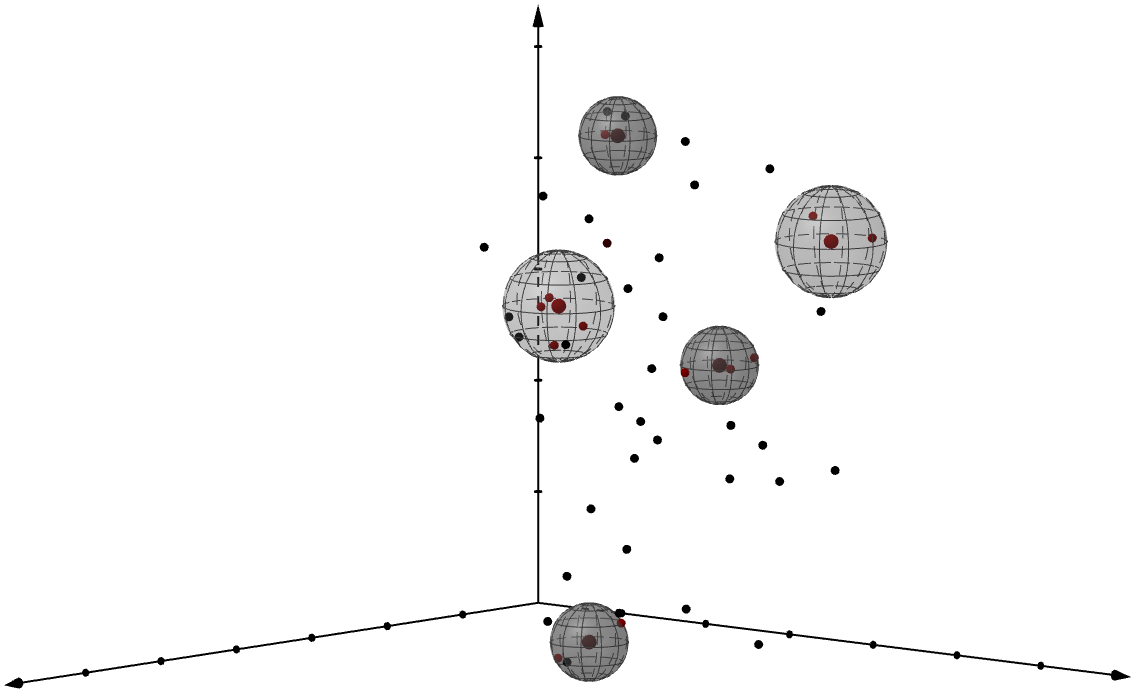}
				\caption{Solution of $\mathbf{p}$-$\rm MTMCLP$ for the 3-dimensional instance of Example \ref{ex:2}.\label{fig1}}
		\end{figure}
		
		In the figure, the covered areas for the discrete facilities are drawn in light gray color while those of the continuous facilities are drawn in dark gray. Red dots indicate covered points and black dots the uncovered ones. The solution achieves an overall demand coverage of $32\%$ being $16\%$ covered by each of the two types of facilities.
\end{example}

\subsection{A `Natural' Non-Linear Model}
 
The modeling framework \eqref{p1...pTMCLP} derived in the previous section is of course valid in the hybridized discrete-continuous setting we are considering. 
The only distinguishing aspect is that in the definition of a set $\mathcal{X}^{(t)} \subseteq \mathbb{S}^{(t)}$ for the finite spaces (i.e., $t=1,\dots,T_1$) we must now consider location-specific radii namely, $\rho_j(t)$ for location $b^{(t)}_j \in \mathbb{S}^{(t)}$, $j \in M^{(t)}$.

For the sake of deriving a suitable mathematical programming formulation for the problem we introduce the following decision variables:

$
y^t_j = 
	\begin{cases} 1, & \text{if facility $b^{(t)}_j \in \mathbb{S}^{(t)}$ is selected}, \\
	0, & \text{otherwise},\end{cases} \quad \text{ for all } t \in \mathcal{T}_1,\: j \in M^{(t)}.
$

$
x^t_i = \begin{cases} 1, & \text{if demand node $a_i$ is covered} \\[-5pt] 
& \text{by the facilities located in $\mathbb{S}^{(t)}$,}\\
0, & \text{otherwise,}\end{cases} \quad \text{for all } t \in \mathcal{T}_1,\: i \in N.
$

$
z^t_{ik} = \begin{cases} 1, & \text{if node $a_i$ is covered} \\[-5pt] 
& \text{by the $k$-th facility in $\mathbb{S}^{(t)}$},\\
0, & \text{otherwise},\end{cases} \quad \text{ for all } i \in N,\: t \in \mathcal{T}_2,\: k \in P_t.
$

$
X^t_k \in\mathbb{S}^{(t)}: \mbox{Coordinates of the $k$-th out of the $p_t$ facilities located in $\mathbb{S}^{(t)}$}, \text{ for } t \in \mathcal{T}_2,\: k \in P_t.$

Using the above decision variables, \eqref{p1...pTMCLP} can be formulated as follows:

\begin{subequations}
	\makeatletter
	\def\@currentlabel{$\mathbf{p}$-{\rm MTMCLP}$^{NL}$}
	\makeatother
	\label{p_1,...,p_TMCLP}
	\renewcommand{\theequation}{$\mathbf{p}$-{\rm MTMCLP}$^{NL}_{\arabic{equation}}$}
	\begin{align}
	\max \quad & \dsum_{i \in N} c_i \left[ \sum_{t\in \mathcal{T}_1}x^t_i + \dsum_{t\in \mathcal{T}_2}\dsum_{k=1}^{p_t} z^t_{ik} \right], \label{P0-omega:obj} \\
	\text{s. t.} \quad & \dsum_{j \in M^{(t)}} y^t_j = p_t,\quad \forall t \in \mathcal{T}_1, \label{P0-omega:pFacilities} \\
	& x^t_i \leq \dsum_{j \in M^{(t)}:\atop \|a_i-b^{(t)}_j\|^{(t)} \leq \rho_j(t)} y^t_j,\quad \forall i \in N,\: \forall t \in \mathcal{T}_1, \label{P0-omega:coveringPresent}\\
	& \dsum_{t=1}^{T_1} x^t_{i} + \dsum_{t=T_1+1}^{T} \dsum_{k=1}^{p_t} z^t_{ik} \leq 1, \quad \forall i \in N, \label{P0-omega:coveringFuture} \\
	& \|a_i- X^t_k\|^{(t)} \leq \rho(t) + U (1 - z^t_{ik}), \label{P0-omega:distance}\\
	&\qquad\qquad\qquad\qquad \forall i \in N,\: \forall t \in\mathcal{T}_2,\:\forall k \in P_t, \nonumber \\
	& x^t_i \in \{0,1\}, \quad \forall t \in \mathcal{T}_1,\: \forall i \in N, \label{P0-omega:domain-x} \\
	& y^t_j \in \{0,1\}, \quad \forall t \in \mathcal{T}_1,\: \forall j \in M^{(t)}, \label{P0-omega:domain-y} \\
	& z^t_{ik}  \in \{0,1\}, \quad \forall i \in N,\: \forall t \in \mathcal{T}_2,\: \forall k \in P_t, \label{P0-omega:domain-z} \\
	& X^t_{k} \in \mathbb{S}^{(t)}, \quad \forall t \in \mathcal{T}_2,\: \forall k \in P_t. \label{P0-omega:domain-X}
	\end{align}
\end{subequations}
In the above model, the objective function \eqref{P0-omega:obj} measures the weighted coverage of the nodes in $\A$ using either the discrete or the continuous facilities; Constraints~\eqref{P0-omega:pFacilities} ensure that exactly $p_t$ discrete facilities are selected in the finite set $\mathbb{S}^{(t)}$; Inequalities~\eqref{P0-omega:coveringPresent} state that node $i$ is covered by a discrete facility iff there is such a facility covering it that has been opened.
Inequalities~\eqref{P0-omega:coveringFuture} guarantee that at most one open facility ``nominated'' for each node and thus each weight $c_i$ is accounted for at most once in the objective function. 
Constraints~\eqref{P0-omega:distance} ensure the proper definition of the $z$-variables. 
In these constraints, $\|a_i-X_k\|^{(t)}$ denotes the $\|\cdot\|^{(t)}$-based distance between point $a_j$ and the $k$-th continuous facility in $\mathbb{S}^{(t)}$ and $U$ is a large enough value. 
Finally, constraints~\eqref{P0-omega:domain-x}--\eqref{P0-omega:domain-X} state the domain of the decision variables.

The large $U$ can be fine tuned.
We can consider a common $U$ for all $t\in \mathcal{T}$. In this case, any value $U \geq \max_{i,j \in N \atop t\in \mathcal{T}} \| a_i - a_j \|^{(t)}$.
Alternatively, for every $t \in \mathcal{T}$ we can set a $t$-specific value $U_t \geq \max_{i,j \in N} \| a_i - a_j \|^{(t)}$.

Observe that for $p_2=\dots,p_T=0$, the problem to solve is a classical discrete $p_1$-MCLP, which can be formulated as above but omitting all terms involving facilities of types $2,\dots,T$ and thus removing all the $z$- and $X$-variables. 
The case $p_t = 0 \ \forall t \in \mathcal{T}\backslash \{T_1+1\}$ and $p_{T_1+1} = 1$ reduces to the classical continuous $p_{T_1+1}$-MCLP which can be also formulated using a reduced set of variables.

\begin{remark}\label{remark:EuclideanDistance}
	\eqref{p_1,...,p_TMCLP} is a Mixed-Integer Non-Linear Programming (MINLP) problem because of \eqref{P0-omega:distance}. 
	In case for some $t\in \mathcal{T}_2$ $\|\cdot\|^{(t)}$ is the Euclidean norm, it is well known that such a type of constraints can be re-written as the following set of linear and second-order cone inequalities:
	
	\begin{align}
	v^t_{ikl} \geq a_{il} - X^t_{kl}, & \: \forall i \in N, \: \forall k \in P_t, \: \forall l \in \{1, \ldots, d\}, \label{P-omega:absolute-1}\\
	v^t_{ikl} \geq -a_{il} + X^t_{kl},& \: \forall i \in N, \: \forall k \in P_t, \: \forall l \in \{1, \ldots, d\},\label{P-omega:absolute-2}\\
	s^t_{ik} \geq \dsum_{l=1}^d \Big(v^t_{ikl}\Big)^2, & \: \forall i \in N, \: \forall k \in P_t,\label{P-omega:euclidean-1}\\
	s^t_{ik} \leq \rho(t) + U\,(1-z^t_{ik}), & \: \forall i \in N, \: \forall k \in P_t, \label{P-omega:distance1}\\
	v^t_{ikl} \geq 0, & \: \forall i \in N, \: \forall k \in P_t, \: \forall l \in \{1, \ldots, d\},\\
	s^t_{ik} \geq 0, & \: \forall i \in N, \: \forall k \in P_t,
	\end{align}
	where $v^t_{ikl}$ represents (via \eqref{P-omega:absolute-1} and \eqref{P-omega:absolute-2}) the absolute value $|a_{il}-X^t_{kl}|$ (here $a_{il}$ and 
	%$X_{il}$ 
	$X^t_{kl}$ stand for the $l$-th coordinate of demand point $a_i$ and the $k$-th facility selected in $\mathbb{S}^{(t)}$, respectively) and $s^t_{ik}$ determines the euclidean distance (by \eqref{P-omega:distance1}) as the squared sum of the absolute differences between the coordinates of $a_i$ and $X^t_{k}$.
	
Consequently, if all the types of continuous facilities use the Euclidean norm, \eqref{p_1,...,p_TMCLP} simplifies to a 
%Mixed Integer Second Order Cone Programming problem, 
mixed-integer second-order cone programming problem, which can be solved using 
%any of the available off-the -shelf softwares.
any available off-the-shelf solver.
\end{remark}

\begin{remark}\label{remark:otherDistances}
	In Remark \ref{remark:EuclideanDistance} one can replace Euclidean distances by block-norm based distances (deriving linear programming models) or by $\ell_\tau$-norm (with $\tau\geq 1$) based distances inducing again mixed-integer second-order cone optimization problems. 
	One may even consider mixed distances (one for each demand point, if one desires to model different coverage areas).
	The interested reader should refer to \citet{BlancoPuertoElHajBenAli:2014} for further details.
\end{remark}

\begin{remark}\label{remark:otherRegions}
Apart from \textit{regular} coverage areas represented by convex surfaces, one could also represent non-convex coverage areas by means of unions of convex (second-order cone representable sets). This type of sets can be efficiently represented using disjunctive constraints that are usually modeled through binary variables \citep[see][for further details on suitable representations of these regions in a location problem]{dolu2020solution}.
\end{remark}

\subsection{An Integer Linear Optimization Model}

The above mixed integer non-linear model becomes intractable for medium or large size instances of the problem even if Euclidean distances are considered. 
Therefore, to successfully tackle the problem other possibilities must to be considered.

Next, we derive an integer linear model based upon projecting out the $X$-variables---which represent the coordinates of the services---by ensuring that these can be \emph{easily} found  (in poly-time) once the different sets of demand points allocated to the same facility are known. 

In what follows, we impose that every selected facility must cover at least one demand point, which, we believe, is a reasonable assumption in practice. 
Furthermore, we consider $T_1=T_2=1$, i.e., a single type of discrete and a single type of continuous facilities are to be located.
The presented results and models can be adapted to the general case. Nevertheless, this would make the contents of this section significantly more involved but not more informative. 
Accordingly, that there are $p_1$ discrete and $p_2$ continuous facilities to locate.

Let us assume that $\mathbb{S}^{(2)}=\R^d$.% is a metric space.
In such a space we define the $\|\cdot\|^{(2)}$-ball centered at $a\in\R^d$ and radius $\rho(2)$ as
$$\B_{\rho(2)}(a) = \{ X \in \R^d \: : \: \|X-a\|^{(2)} \leq \rho(2) \}.$$

The new model we propose relies on Lemma~\ref{lem:1} whose proof is straightforward:
\begin{lem}\label{lem:1}
	Let $N_1,\dots,N_{p_2} \subseteq N$ be $p_2$ nonempty disjoint subsets of $N$. 
	Then, they induce a solution to the $(p_1,p_2)-{\rm MTMCLP}$ in the $X$-variables (where $N_k$ are the points covered by $X_k$) if and only if
	$$
	\bigcap_{i \in N_k} \B_{\rho(2)}(a_i) \neq \emptyset, \forall k \in P_2.
	$$
\end{lem}
Lemma~\ref{lem:1} allows us to rewrite constraints \eqref{P0-omega:distance} (for $T=2$) as linear constraints and thus to formulate the \emph{Hybridized  $(p_1,p_2)$-Maximal Covering Location problem} (HMCLP, for short) as follows:
\begin{subequations}
	\makeatletter
	\def\@currentlabel{${\rm HMCLP}^{IP}$}
	\makeatother
	\label{prMCLP-H}
	\renewcommand{\theequation}{${\rm HMCLP}^{IP}_{\arabic{equation}}$}
	\begin{align}
	\max \quad & \dsum_{i \in N} c_i \left[ x_i + \dsum_{k \in K} z_{ik} \right], \label{P0-omega:objH} \\
	\text{s. t.} \quad & \dsum_{j \in M} y_j = p_1, \label{P0-omega:pFacilitiesH} \\
	& x_i \leq \dsum_{j \in M \atop \|a_i-b_j\|^{(1)} \leq \rho_j(1)} y_{j},\quad \forall i \in N, \label{P0-omega:coveringPresentH}\\
	& x_{i} + \dsum_{k=1}^{p_2} z_{ik} \leq 1, \quad \forall i \in N, \label{P0-omega:coveringFutureH} \\
	& \dsum_{i \in Q} z_{ik} \leq |Q|-1, \forall k \in P_2,\: \forall Q \subseteq N: \bigcap_{i \in Q} \B_{\rho(2)}(a_i) = \emptyset, \label{P0-omega:distanceH} \\
	& x_{i} \in \{0,1\}, \quad \forall i \in N, \nonumber\\
	& y_{j} \in \{0,1\}, \quad \forall j \in M, \nonumber\\
	& z_{ik}  \in \{0,1\}, \quad \forall i \in N,\: \forall k \in P_2.  \nonumber
	\end{align}
\end{subequations}
In this model, we have simplified some notation namely by removing the type index from the decision variables as well as from set $M$ since we only have one type of discrete and one type of continuous facilities. Also in the above model, \eqref{P0-omega:distanceH} ensure that the set of points covered by a continuous facility verifies the condition of Lemma~\ref{lem:1},  i.e., sets of \textit{incompatible} demand points are not allowed to be allocated to the same continuous facility. This constraint replaces the non-linear constraint 
in \eqref{p_1,...,p_TMCLP}.

Moreover, the above model does no include the variables $X_k$. 
In fact, once an optimal solution is obtained for \eqref{prMCLP-H}, we can use the values of the $z$-variables, say $\{\bar{\mathbf{z}}\}$, to find explicit optimal coordinates for the new facilities to install. In particular, the coordinates of the $k$-th facility to install ($k\in P_2$) can be given by any vector $X_k$ satisfying:
$$
X_k  \in \bigcap_{i \in N : \atop \bar z_{ik}=1} \B_{\rho(2)}(a_i), \text{for all $k \in P_2$. }
$$

A center $X_k$ in the above intersection can be found in polynomial time either solving a second order cone optimization problem or by solving a one-center facility location problem.

Solving \eqref{prMCLP-H} requires incorporating exponentially many constraints---inequalities \eqref{P0-omega:distanceH}.
Interestingly, the above formulation can be simplified (reducing from exponential to polynomially many constraints in the form of \eqref{P0-omega:distanceH}) by means of Helly's Theorem \citep{Helly:1923} \citep[see also][]{DanzerGruenbaumKlee:1963}.
Invoking that result, provided that the continuous space is $\R^d$, only intersections of $(d+1)$-wise balls are needed to check:
$$
\B_{\rho(2)}(a_{i_1}) \cap \dots \cap \B_{\rho(2)}(a_{i_{d+1}}),
$$
for all $a_{i_{1}}, \ldots, a_{i_{d+1}} \in \mathcal{A}$.

Despite this simplification, the number of constraints may still be large and thus making the problem more difficult to solve. 
Instead, the problem can be tackled by considering an incomplete formulation (removing \eqref{P0-omega:distanceH}) and iteratively incorporating these constraints \textit{on-the-fly}, as needed.

The selection of the constraints to incorporate in each iteration is found using the following separation strategy: 
After solving \eqref{prMCLP-H} with none or part of the constraints \eqref{P0-omega:distanceH} a solution, say $\bar{\mathbf{z}}$, is obtained.
Then, for each $k \in P_2$ the define set $Q_k= \{i \in N: \bar z_{ik} =1\}$. 
One can check for the validity of the set $Q_k$ as a feasible cluster of demand points for our problem by solving the $1$-center problem for the points in such a set.
In case the optimal coverage radius obtained is less than or equal to $\rho(2)$, one knows that $Q_k$ is a valid subset of demand points that can be covered by the same server. 
Otherwise, the solution violates the relaxed constraints, and thus we add the cut
\begin{equation}\label{Qk}
\dsum_{i \in Q_k} z_{ik'} \leq |Q_k|-1, \forall k' \in P_2,
\end{equation}
to ensure that such a solution is no further deemed feasible and thus obtained again. 

The $1$-center problem with Euclidean distances on the plane is known to be solvable in polynomial time 
%using classical continuous location algorithms 
\citep[see e.g.,][]{EH72}. Extensions to higher dimensional spaces and generalized covering shapes have been recently proved to be also poly-time solvable \citep{BP21}.

The above procedure can be embedded into a branch-and-cut approach by means of lazy constraints.

The following result holds, which helps finding (and thus ignoring) dominated cuts:

\begin{prop}\label{prop:EH}
	Let $\mathbf{z} \in \{0,1\}^{n\times p_2}$ and $Q, Q' \subseteq N$ be such that $Q \subset Q'$. Then, if $\mathbf{z}$ violates \eqref{P0-omega:distanceH} for the set $Q$ violates, then, $\mathbf{z}$ also violates the constraint forthe set $Q'$. Thus, the cut induced by $Q$ strictly dominates the one induced by $Q'$.
\end{prop}
\begin{proof}
	We suppose that $Q$ violates the constraint \eqref{P0-omega:distanceH}, this means $\bigcap_{i \in Q} \B_{\rho(2)}(a_i) = \emptyset$ and we get $\dsum_{i \in Q} z_{ik}  > |Q|-1$. Therefore,
	$$
	\dsum_{i \in Q'} z_{ik}  =  \dsum_{i \in Q} z_{ik} +  \dsum_{i \in Q'\backslash Q} z_{ik} > |Q|-1 + (|Q'|-|Q|) = |Q'|-1, \forall k \in P_2.
	$$
	
	Then, we have that $Q'$ also violates the constraints and the cut induced by $Q$ strictly dominates the one induced by $Q'$.
\end{proof}

In the next section we provide further details on how we can take advantage from these contents using a more specific setting.

\section{The Particular Case of the Euclidean Plane}
\label{sec:EuclideanPlane}

In this section we focus on the particular case in which the continuous facilities are to be located in the Euclidean plane. 
This allows deepening the discussion already presented and also to consider an alternative model for the problem.

\subsection{A Branch-and-Cut Algorithm based on ${\rm (HMCLP}^{IP}${\rm )}}
\label{subsec:BandC}
Let us take again the integer linear model (${\rm HMCLP}^{IP}$).
In the case of the plane ($d=2$), the application of Helly's Theorem described above guarantees that we only need to check intersection or $3$-wise balls. In particular, we need to check if such intersections are empty. 
If so, we incorporate the adequate constraints to avoid searching for facilities in those intersections. Note that is true for every norm we adopt.
Two cases may emerge for any set of three demand points $\{a_{i_1}, a_{i_2}, a_{i_3}\}  \subset \A$:
\begin{description}
	\item[Case 1:] $\B_{\rho(2)}(a_{l_1}) \cap \B_{\rho(2)}(a_{l_2}) = \emptyset$ for some $l_1, l_2 \in \{i_1, i_2, i_3\}$.
	
	In this case, points $a_{l_1}$ and $a_{l_2}$ are incompatible---they cannot be covered by the same center.
	Hence impose
	\begin{equation}\label{Helly-1}\tag{2-Wise}
	z_{{l_1}k} + z_{{l_2}k} \leq 1, \forall k \in P_2.
	\end{equation}
	
	\item[Case 2:] The pairwise intersections are non-empty but $\B_{\rho(2)}(a_{i_1}) \cap \B_{\rho(2)}(a_{i_2}) \cap \B_{\rho(2)}(a_{i_3})  = \emptyset$.
	
	In this case, the three points cannot be covered by the same facility and thus we impose
	\begin{equation}\label{Helly-2}\tag{3-Wise}
	z_{{i_1}k} + z_{{i_2}k} + z_{{i_3}k} \leq 2, \forall k \in P_2.
	\end{equation}
\end{description}

Constraints \eqref{Helly-1} and \eqref{Helly-2}  for all subsets of three points in $\A$ replace the constraint \eqref{P0-omega:distanceH} in formulation \eqref{prMCLP-H} in the planar case. Although these constraints (in worst case) are $O(n^3)$, most of them are needless to construct the optimal solution of the problem. 

The problem can now be solved considering an incomplete formulation (removing constraints \eqref{Helly-2} from \eqref{prMCLP-H}) and incorporating these constraints \textit{on-the-fly}, as needed using a similar strategy than the one used for the general formulation and described above. It avoids checking the three-wise intersections of points which can be computationally costly for large instances.
We embed this approach into a Branch-and-Cut scheme which is reinforced by the following  elements:
\begin{description}[style=unboxed,leftmargin=0cm]
	\item[Separation Oracle] Given an optimal solution $\bar{\mathbf{z}}$ to the incomplete model, for every $k \in P_2$ we consider the set $Q_k= \{i \in N: \bar z_{ik} =1\}$. 
	The validity of the cluster $Q_k$ is checked using the algorithm proposed by \citet{EH72} that computes (in polynomial time), the center and the minimum radius covering the demand points in $Q_k$.
	
	To make this work self-contained, we recall that the algorithm proposed by \citet{EH72} is based on the construction of  disks covering three points until the whole set is covered.
	Thus is accomplished sequentially in such a way that the covering radius increases at each step of the procedure. 
	
	We take further advantage from this procedure to incorporate more than a single cut in each iteration of the branch-and-cut algorithm.
	Considering a set of points $Q_k$, if at some step when applying the algorithm by \citet{EH72} the coverage radius becomes larger than $\rho(2)$ (it is not feasible to clustering the points in $Q_k$ to be served by the same server), then all the sets of two or three points used to construct the minimum enclosing disks so far are used to generate constraints \eqref{Helly-1} and \eqref{Helly-2}.
	This strategy has proven to alleviate the resolution of the problem.
	
	\item[Initial Pool of Constraints] 
	%Apart from the constraints defining the incomplete formulation that we consider, we have designed a procedure to generate an initial pool of constraints in the form of \eqref{Helly-2}. 
	We propose procedure to generate an initial pool of constraints in the form of \eqref{Helly-2} to be included in the initial incomplete model.
	First, we construct clusters of demand points with maximum distance between two clustered points fixed to $\rho(2)+\varepsilon$, with $\varepsilon > 0$ (implemented in the Python module \texttt{scipy} through the function \texttt{fcluster}). At each of these clusters, we checked the validity of them as a solution of our MCLP using the same strategy used in our separation oracle. In case some of them is not valid,  we incorporate the pool all the constraints of type  \eqref{Helly-2} that are violated. The procedure is repeated for different values of $\varepsilon$.
	\item[Symmetries] Our formulation \eqref{prMCLP-H} is highly symmetric in the sense that any permutation of the $k$ indices in the $z$-variables results in an alternative solution (with the same objective value). 
	To break symmetries in our model thus hoping to speed up the resolution of the problem, we incorporate the following constraints that we can observe to be straightforwardly valid for the problem:
	\begin{equation}
	\dsum_{i \in N} w_i z_{i (k-1)} \leq \dsum_{i \in N} w_i z_{ik}, \ \forall k \in P_2\backslash\{1\}.
	\end{equation}
	With these constraints, among all the posible sorting of facilities, we choose one producing a non-decreasing sequence of covered demand.
\end{description}

\subsection{An Alternative IP Model}
\label{subsec:alternativeIP}
As largely explained, model (${\rm HMCLP}^{IP}$) is quite general;
it is valid in metric spaces of any dimension $d \geq 2$ and for every distance of interest.
In $2$-dimensional spaces and in addition to the developments presented in Section  \ref{subsec:BandC}, we can derive an alternative Integer Linear Programming formulation.
This is accomplished by finding a finite dominating set, which is done using the discretization technique proposed by \citet{Church:1984} for the MCLP.

Let $\mathcal{B}$ be the set consisting of the demand points in $\A$ and also the intersection points of the pairwise intersections of the boundary of the $\|\cdot\|^{(2)}$-balls centered at the demand points with radius $\rho(2)$, i.e.,
$$
\mathcal{B} = \A \cup \bigcup_{i, l\in N:\atop i< l} \Big(\partial \B_{\rho(2)} (a_i) \cap \partial \B_{\rho(2)} (a_l)\Big),
$$
where $\partial \B_{\rho(2)}(a)$ stands for the boundary of the ball $\B_{\rho(2)}(a)$. 
Inspired by the terminology used by \citet{Church:1984} the set $\mathcal{B}$ will be designated by the Balls Intersection Points Set (BIPS).

\begin{lem}
	There exists an optimal solution to $(p_1,p_2)$-MTMCLP where the continuous facilities belong to the set $\mathcal{B}$.
\end{lem}
\begin{proof}
	Let $\mathcal{X}^2 = \{\overline{X}_1, \ldots, \overline{X}_{p_2}\}$ be optimal positions for the continuous facilities. Clearly, $\overline{X}_k$ belongs to the intersection of the balls centered at the covered points, i.e.,$$\overline{X}_k \in \bigcap_{i \in N \: : \: \|a_i-\overline{X}_k\|^{(2)} \leq \rho(2)} \B_{\rho(2)}(a_i).$$ 
	Thus, we can replace $\overline{X}_k$ by any of those intersection points keeping the same coverage level.
	Clearly, these points belong to $\mathcal{B}$.
\end{proof}

Considering the set $\mathcal{B}$,  it is possible to reformulate  $(p_1,p_2)$-MTMCLP$^{NL}$ as an integer linear programming problem, in which the selection of the continuous facilities is replaced by the search of the optimal $\mathcal{B}$-points to open. 
Let us denote by $\mathcal{B}=\{\gamma_1, \ldots, \gamma_{|\mathcal{B}|}\}$ and $L=\{1, \ldots, |\mathcal{B}|\}$, and consider the following sets of decision variables:

\medskip
$y^1_{j} = \begin{cases} 1, & \text{if facility $b_j$ is selected}, \\ 0, & \text{otherwise},\end{cases}
$ $\forall j \in M$,

\medskip
$y^2_{l} = \begin{cases} 1, & \text{if point $\gamma_l \in \mathcal{B}$ is selected}, \\ 0, & \text{otherwise},\end{cases} \quad \forall l \in L.
$

\medskip
The problem can be reformulated as follows:
\begin{subequations}
	\makeatletter
	\def\@currentlabel{${\rm HMCLP}^{BIPS}$}
	\makeatother
	\label{p1p2MCLP_CIP}
	\renewcommand{\theequation}{${\rm HMCLP}^{BIPS}_{\arabic{equation}}$}
	\begin{align}
	\max \quad & \dsum_{i \in N} c_ix_i \label{P0-omega:objCIP} \\
	\text{s. t.} \quad & \dsum_{j \in M} y^1_j = p_1, \label{P0-omega:pFacilitiesCIPa} \\
	& \dsum_{l \in L} y_l^2 = p_2, \label{P0-omega:pFacilitiesCIPb} \\
	& x_i \leq \dsum_{j \in M \atop \|a_i-b_j\|^{(1)} \leq \rho_j(1)} y^1_{j} +  \dsum_{l \in L \atop \|a_i-c_l\|^{(2)} \leq \rho(2)} y^2_{l},\quad \forall i \in N, \label{P0-omega:coveringPresentCIP}\\
	& x_{i} \in \{0,1\}, \quad \forall i \in N, \nonumber \\%\label{P0-omega:domain-xCIP}
	& y^1_{j} \in \{0,1\}, \quad \forall j \in M, \nonumber \\%\label{P0-omega:domain-yCIPa} \\
	& y^2_{l} \in \{0,1\}, \quad \forall l \in L, \nonumber%\label{P0-omega:domain-yCIPb}
	\end{align}
\end{subequations}
where constraint \eqref{P0-omega:pFacilitiesCIPa} enforces opening exactly $p_1$ of the facilities from $\mathbb{S}^{(1)}$, while constraint \eqref{P0-omega:pFacilitiesCIPb} ensures opening exactly $p_2$ of the additional (continuous) facilities in $\mathbb{S}^{(2)}=\R^2$. 
Constraints \eqref{P0-omega:coveringPresentCIP} allow to determine whether a demand point is covered or not by any of the available open facilities (from $\mathbb{S}^{(1)}$ or $\mathcal{B}$). 
Note that in the above model the domain of the $x$-variables can be relaxed to the interval $[0,1]$.

The problem above is a particular version of the classical Discrete Maximal Covering Location problem in which two different types of facilities are desired to be open, $p_1$ of type discrete and $p_2$ of type continuous. 

%One major difficulty underlying the above model stems from the dimension of the underlying continuous space. 
A major issue of concern in the above model is the number of $y^2$-variables, which coincides with the number of points in $\mathcal{B}$.
This number is of order $O(n^2)$ considering one type of continuous facilities but will be of course larger if additional continuous facilty types exist since we must add additional sets of $y-$variables---one for each facility type.

Still considering a single type of continuous facilities, the size of the set $\mathcal{B}$ can be reduced following the strategy proposed in \citet{Church:1984}.
However, again we must note that that author worked only with $\ell_1$- and $\ell_2$-norm for which a dominance relation between the points allows removing some elements from $\mathcal{B}$.
In our case, although working in the Euclidean plane, we can consider distances other than the $\ell_1$-norm and the $\ell_2$-norm.
This, again, creates some challenges to the above model since finding the BIPS is far from straightforward. 

Overall, the model just proposed can be promising under a particular case: only one type of continuous facilities exist and the distances of interest reduce to $\ell_1$ or $\ell_2$ norms.

\section{Computational Experiments}
\label{sec:ComputationalExperiments}

In this section we report on the results of a series of computational experiments performed to empirically assess our methodological contribution for the hybridized MCLP presented in the previous sections.

\subsection{The Test Data}
To run the experiments, we made use of real geographic and demographic information from Manhattan Island, NY, USA.
The data was collected from \href{https://data.cityofnewyork.us/Housing-Development/Housing-New-York-Units-by-Building/hg8x-zxpr}{\url{data.cityofnewyork.us}}. 
The main instance consists of the (planar) geographical coordinates of the 920 main buildings on the island with demand weights given by the number of people living in at each of the buildings. 
The complete data set used in our experiments is available in our repository \url{github.com/vblancoOR/MTMCLP}.

\begin{figure}[H]
	\begin{subfigure}[b]{.25\linewidth}
		\centering
		\includegraphics[scale=0.7]{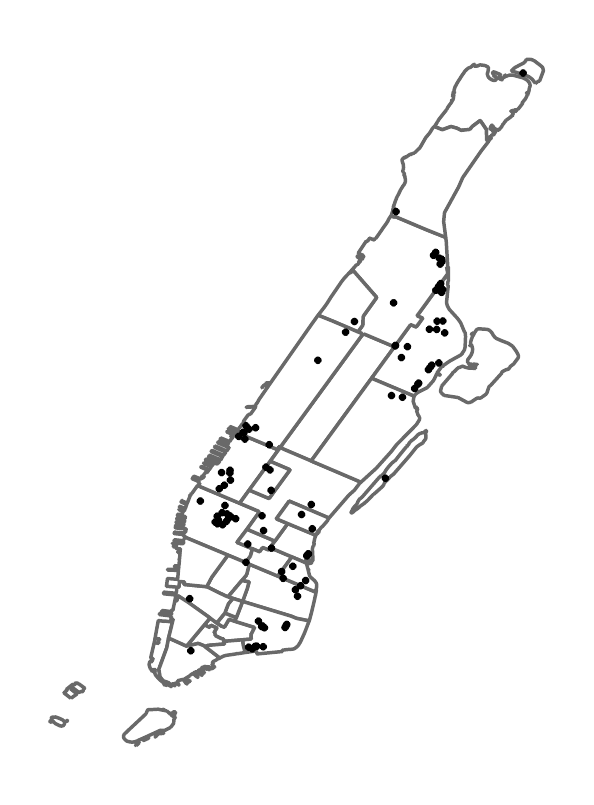}
		\caption{$n = 100$.}
	\end{subfigure}~\begin{subfigure}[b]{.25\linewidth}
		\centering
		\includegraphics[scale=0.7]{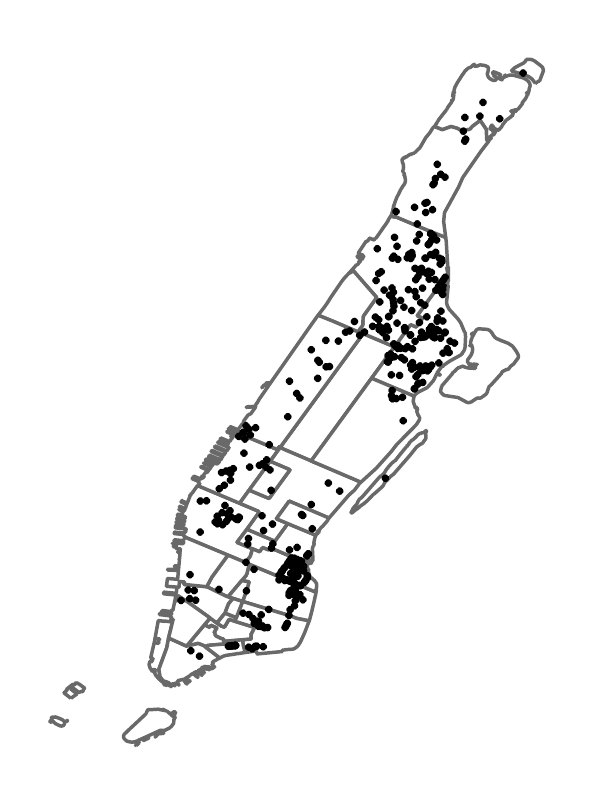}
		\caption{$n = 400$.}
	\end{subfigure}~\begin{subfigure}[b]{.25\linewidth}
		\centering
		\includegraphics[scale=0.7]{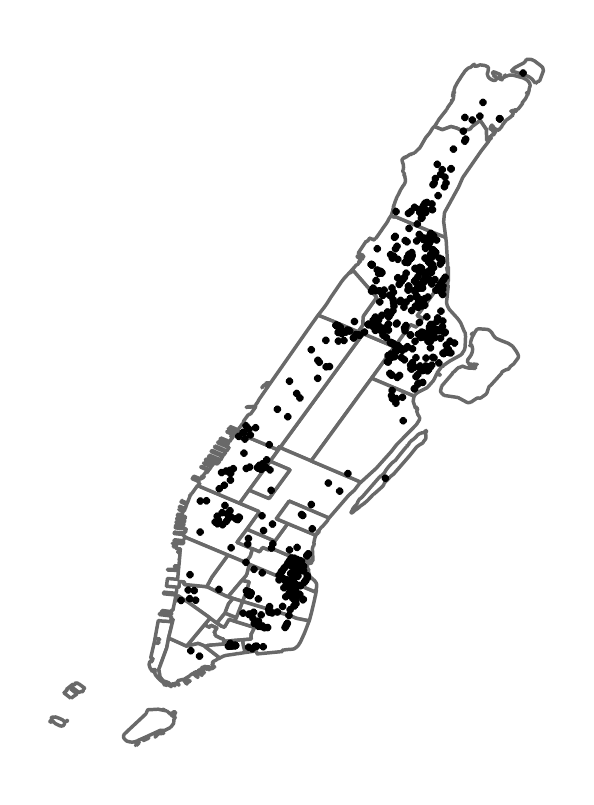}
		\caption{$n = 700$.}
	\end{subfigure}\begin{subfigure}[b]{.25\linewidth}
		\centering
		\includegraphics[scale=0.7]{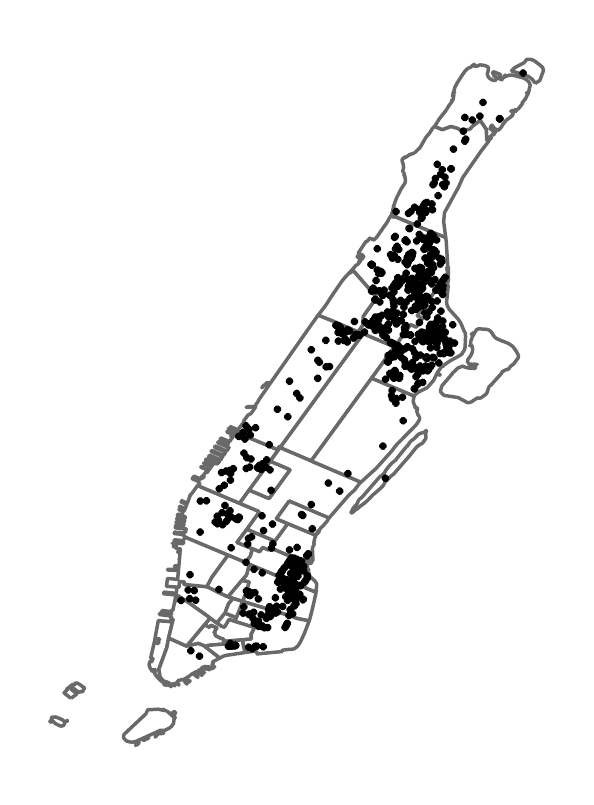}
		\caption{$n = 920$.}
	\end{subfigure}
	\caption{Different instances considered from the Manhattan dataset.}
	\label{f:manhattan}
\end{figure}

To test the scalability of the problem we are investigating, different subsets of the complete data we considered each with a different size. 
We sorted the location indices of the buildings according to demands and we made subsets from the buildings with the largest demands. 
We have considered subsets of cardinality $n \in \{400,500,700,920\}$ with 920 corresponding to the complete data set.
Figure \ref{f:manhattan} depicts the different sizes considered for the Manhattan data set.

We assume one type of discrete facility and one type of continuous facility.
Additionally, we suppose that the potential facility location for the discrete facilities correspond to buildings underlying the instance of the problem and we consider that the coverage radius is equal for all of them. 
As for facilities to locate we adopted $p_1, p_2 \in \{1,2,3,4,5\}$.
Finally, we consider two radii for each type of facilities: for discrete facilities we consider $\rho(1) \in \{0.008,0.012\}$ and for the continuous radii the values are $\rho(2) \in \{0.005,0.01\}$. These radii are adjusted to the real coordinates (latitude and longitude) of the demand points and are equivalent to $\{810,1280\}$ (for $\rho(1)$) and $\{430,1050\}$ (for $\rho(2)$) meters.
Combining all the parameters, a total of 400 instances have been tested (four values for $n$, five values for $p_1$, five values for $p_2$, two values for $\rho_1$ and two values for $\rho_2$).

In the tests whose results we detail next, we run: (i) the non-linear formulation \eqref{p_1,...,p_TMCLP}, (ii) the integer programming formulation \eqref{prMCLP-H} (with constraints \eqref{Helly-1} and \eqref{Helly-2} and the Branch-and-Cut approach derived from an incomplete formulation of \eqref{prMCLP-H} in which only the pair-wise intersection constraints are considered, and (iii) the integer programming formulation \eqref{p1p2MCLP_CIP}, 

All the experiments have been run on a virtual machine in a physical server equipped with 12 threads from a processor AMD EPYC 7402P 24-Core Processor, 64 Gb of RAM and running a 64-bit Linux operating system. 
The models were coded in Python 3.7 and we used Gurobi 9.1 as optimization solver. A time limit of 1 hour was fixed for all the instances.

\subsection{Results}
\label{subsec:results}

\begin{figure}[ht]
\captionsetup[subfigure]{aboveskip=-3pt,belowskip=-3pt}
	\begin{subfigure}[b]{.5\textwidth}
		\includegraphics{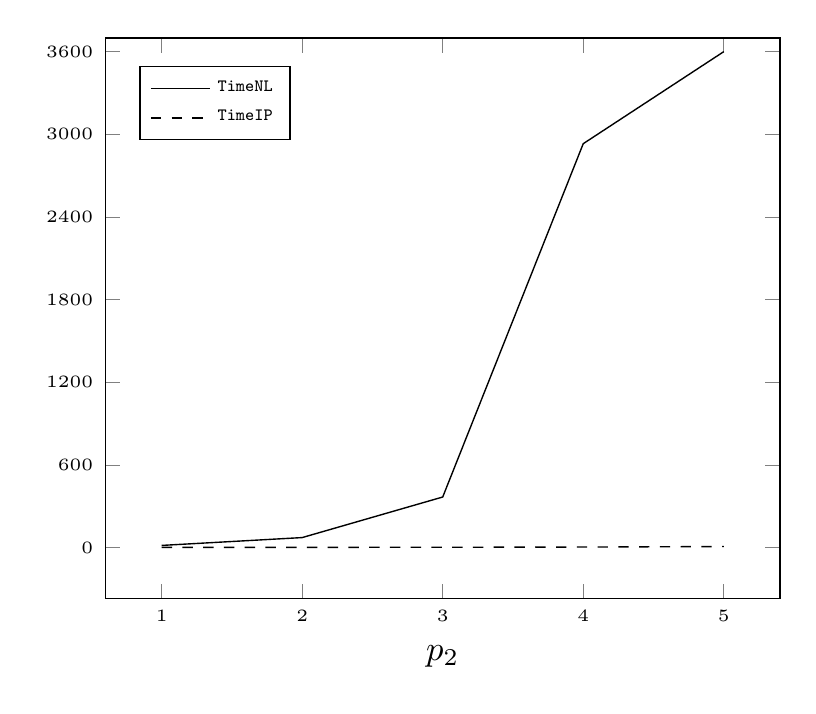}
		\caption{$\rho(2) = 0.005$}
	\end{subfigure}
	\begin{subfigure}[b]{.5\textwidth}
		\includegraphics{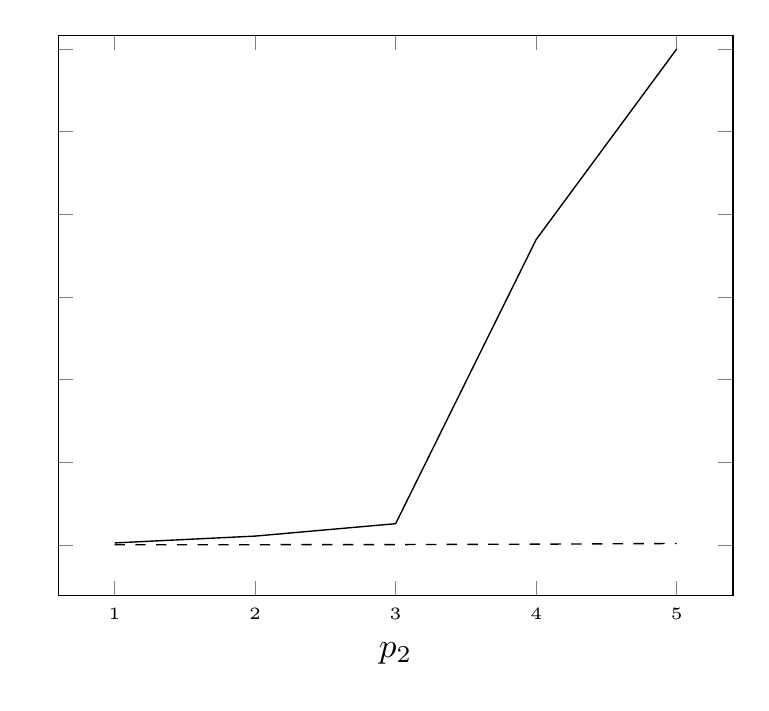}
		\caption{$\rho(2) = 0.01$}
	\end{subfigure}
	\caption{Average of CPU times required by models \eqref{p_1,...,p_TMCLP} (straight line) and \eqref{prMCLP-H} (dashed line) when solving the instances with 100 demand nodes.}
	\label{f:zvsnl}
\end{figure}

Figure~\ref{f:zvsnl} depicts results for the smallest instances with the purpose of comparing the non-linear formulation \eqref{p_1,...,p_TMCLP} and the complete pure integer formulation \eqref{prMCLP-H}. We show in that figure the average CPU times for both approaches by aggregating the different values of $p_1$ and $\rho(1)$. The instances with $n=100$ demand nodes were considered in these tests.

We show how the computational times increase when the number of continuous facilities to locate increases for the two values of $\rho(2)$.
In particular we observe that the CPU time reaches the time limit for $p_2=5$. 
These results give strong evidence to our intuition: the non-linear formulation is computationally demanding even for small instances of the problem.

%\FSG{Francisco: I don not clearly understand the exact values taken for $p_1$. Above we say that $p_1 \in \{1,2,3,4,5 \}$ but then are the results aggregated somehow in Figure \eqref{f:zvsnl}? The comment is also valid for Table~\ref{t:ztimes}. I thought we had worked only with $p_1=5$ as I wrote above. We need to check this.}

Given the results presented in Figure~\ref{f:zvsnl} we decided to proceed the tests without the non-linear model.

Focusing on model \eqref{prMCLP-H}, Table \ref{t:ztimes} contains the number constraints \eqref{Helly-1} and \eqref{Helly-2} generated as well as the CPU time required for that generation.
In this table we observe that the linear model grows a lot when the number of demand points increases, which makes is clearly more difficult to tackle.
Moreover, the CPU time required to check all intersections increases significantly for constraints \eqref{Helly-2}.
Summing up, we easily conclude that making use of a branch-and-cut approach such as the one we propose in Section~\ref{subsec:BandC} is totally advisable in this case since we can avoid the clear burden that corresponds to computing and using all the constraints.
\begin{table}[ht]
\begin{center}
	\adjustbox{max width=0.75\textwidth}{
		\begin{tabular}{lrrrrr}
		\toprule
		             & & \multicolumn{2}{c}{\eqref{Helly-1}} & \multicolumn{2}{c}{\eqref{Helly-2}} \\
		\texttt{n} & \multicolumn{1}{l}{$\rho(\texttt{2})$} & \multicolumn{1}{c}{\texttt{CPUTime}} & \multicolumn{1}{c}{\texttt{\# Inequalities}}  & \multicolumn{1}{c}{\texttt{CPUTime}} & \multicolumn{1}{c}{\texttt{\# Inequalities}}  \\ 
		\midrule
		\multirow{2}{*}{400} & 0.005 & 0.62 & 71780 & 490.70 & 1457 \\ 
		 & 0.01 & 0.55 & 61532 & 2658.11 & 5925 \\ 
		 \midrule
		\multirow{2}{*}{500} & 0.005 & 0.91 & 112404 & 1483.38 & 3701 \\ 
		 & 0.01 & 0.91 & 94279 & 9273.37 & 16600 \\ 
		 \midrule
		\multirow{2}{*}{600} & 0.005 & 1.35 & 161160 & 3997.49 & 6937 \\
		 & 0.01 & 1.29 & 133860 & 25277.75 & 30218 \\ 
		 \midrule
		\multirow{2}{*}{700} & 0.005 & 1.85 & 218260 & 9596.85 & 12171 \\ 
		 & 0.01 & 1.80 & 179211 & 61696.80 & 58717 \\ 
		 \midrule
		\multirow{2}{*}{920} & 0.005 & 3.14 & 369085 & 52788.62 & 40047 \\ 
		 & 0.01 & 3.09 & 290690 & 358966.03 & 184221 \\ 
		 \bottomrule
	\end{tabular}
	}
	\caption{Averaged number of Constraints~\eqref{Helly-1} and \eqref{Helly-2} generated and the CPU time (seconds) required for their generation.}
	\label{t:ztimes}
	\end{center}
\end{table}

Next, we focus on the results corresponding to running the B\&C procedure devised for tackling model \eqref{prMCLP-H} and also the when formulation \eqref{p1p2MCLP_CIP} is adopted.
The common ground for this comparison is the Euclidean plane using the Eucliden distance.
In fact, this is the setting in which we can use model \eqref{p1p2MCLP_CIP}. 
Another possibility would be to use the $\ell_1$-norm but then we would clearly be favor model \eqref{p1p2MCLP_CIP}, which is already done using the Euclidean distance.

In tables \eqref{t:computational1} and  \eqref{t:computational2}  we present the results obtained. Each table refer to a different covering radius for the continuous facilities ($\rho(2) = 0.005$ and $\rho(2)=0.01$, respectively) and are similarly organized. Both tables show the information aggregating the different values of the discrete facilities, $p_1$, and the radii of them, $\rho(1)$. The first three columns provide the details of the instance being solved: the radii used for the continuous facilities, the number of demand nodes, and the number of continuous facilities to locate.
Next, seven blocks of columns are presented.
 In the blocks, the columns contains results for model \eqref{p1p2MCLP_CIP} (\texttt{BIPS}) and the branch-and-cut procedure devised for tackling model \eqref{prMCLP-H} (\texttt{B\&C}). Some blocks only contain information of the branch-and-cut approach since it is not applicable to the \texttt{BIPS} approach.
The first block (columns 4 and 5) gives the overall CPU time in seconds (averaged only on the instances that were solved up to optimality within the time limit) required to solve the problem. This time includes: 
the time required by Gurobi to solve the IP models (second block)  of solved instances; the pre-processing time---time for generating the \texttt{BIPS} in the case of Constraints \eqref{P0-omega:coveringPresentCIP} and the time for building the initial pool of Constraints \eqref{Helly-1} in the case of model \eqref{prMCLP-H} (third block);
the time to generate the difficult constraint \eqref{P0-omega:coveringPresentCIP} and its respective in the \eqref{prMCLP-H} formulation (fourth block);
callback times for the branch-and-cut approach (fifth block); 
percentage gap at termination and number of unsolved instances within the time limit-- out of 10 (sixth  block); and total number of constraints used (sixth block). The percentage gaps and the number of unsolved instances are not reported for the \texttt{BIPS} approach since all the instances were solved up to optimality within the time limit.

In view of the results of tables \eqref{t:computational1} and  \eqref{t:computational2} we draw several conclusions.
First, we realize that the difficulty in using model ${\rm HMCLP}^{BIPS}$ stems from loading Constraints~\eqref{P0-omega:coveringPresentCIP}.
In fact, to accomplish this we need to check all the intersections.
This takes a long time if the size of \texttt{BIPS} is large.
Another challenge for that model regards the time to find the \texttt{BIPS}.
When the number of demand nodes increases and the radius is large enough to require checking all (or a large majority of) points, then it takes a long time.
From the above observations we conclude that in applications with a huge number of points and large radii (e.g. clustering problems easily lead to such cases) model ${\rm HMCLP}^{BIPS}$ may become intractable.
Still concerning this model, we note that the CPU time required to solve it seems quite indifferent to the number of continuous facilities to locate.
 
We turn now our attention to the \texttt{B\&C} algorithm devised for model \eqref{prMCLP-H}. 
From tables \eqref{t:computational1} and  \eqref{t:computational2}  we see that the model becomes more challenging when the number of continuous facilities to locate increase.
Nevertheless, in those cases in which a proven optimal solution could not be found within the time limit imposed, the final gap is quite small.
Still concerning the \texttt{B\&C} procedure we observe that the pre-processing time as well as that for generating violated cuts and for the callback checks are all very low.
This makes the whole algorithm more efficient.
The number of continuous facilities to locate is clearly the factor influencing the most the performance of the algorithm since the number of incorporate cuts increases significantly.
Still, we see that the \texttt{B\&C} algorithm outperforms the plain use of model ${\rm HMCLP}^{BIPS}$ when the number of continuous facilities to locate is small (one or two). 
Nevertheless, we must recall that the comparison that can be observed in both tables can be made only because we are working on the Euclidean plane and using the Euclidean distance.
As discussed above, if this was not the case, then we would certainly need to resort to the \texttt{B\&C} algorithm devised since it is quite insensitive to the adopted norm, which is far from the case when we need to determine \texttt{BIPS}. The results show that the \texttt{B\&C} approach is certainly a viable algorithm, which is quite general.

\begin{table}[htbp]
	\centering
	\adjustbox{width=\textwidth}{\begin{tabular}{ll|rrrrrrrrr|rc|rr}
		\multicolumn{2}{l|}{} &\multicolumn{9}{c|}{\texttt{CPUTime (secs.)}} &\multicolumn{2}{c|}{} &\multicolumn{2}{c}{}\\
			
		\multicolumn{2}{l|}{} & \multicolumn{ 2}{c}{\texttt{Total}} & \multicolumn{ 2}{c}{\texttt{Solving}} & \multicolumn{ 2}{c}{\texttt{Prepr.}} & \multicolumn{ 2}{c}{\texttt{Ctrs. Gen.}} & \multicolumn{ 1}{c|}{\texttt{Callback}} & \multicolumn{ 1}{c}{\texttt{MIPGAP}} & \multicolumn{ 1}{c|}{\texttt{\#Unsolved}} & \multicolumn{ 2}{c}{\texttt{\#Ctrs}} \\  \cline{ 3- 15}
		
		n  & \multicolumn{1}{c|}{$p_2$} & \multicolumn{1}{c}{ \texttt{BIPS} } & \multicolumn{1}{c}{ \texttt{B\&C} } & \multicolumn{1}{c}{ \texttt{BIPS} } & \multicolumn{1}{c}{ \texttt{B\&C} } & \multicolumn{1}{c}{ \texttt{BIPS} } & \multicolumn{1}{c}{ \texttt{B\&C} } & \multicolumn{1}{c}{ \texttt{BIPS} } & \multicolumn{1}{c}{ \texttt{B\&C} } & \multicolumn{1}{c|}{ \texttt{B\&C} } & \multicolumn{1}{c}{ \texttt{B\&C} } & \multicolumn{1}{c|}{ \texttt{B\&C} } & \multicolumn{1}{c}{\texttt{BIPS}} & \multicolumn{1}{c}{ \texttt{B\&C} } \\ \hline
		
 \multirow{5}{*}{400} & 1 & 53.45 & 6.24 & 0.27 & 3.75 & 7.10 & 1.38 & 46.09 & 1.10 & 0 & 0 & 0 & 402 & 72587 \\ 
& 2 & 53.28 & 51.91 & 0.23 & 49.41 & 7.10 & 1.36 & 45.96 & 1.10 & 0.04 & 0 & 0 & 402 & 144373 \\ 
& 3 & 53.46 & 229.64 & 0.24 & 227.00 & 7.10 & 1.37 & 46.12 & 1.12 & 0.15 & 0 & 0 & 402 & 216159 \\ 
& 4 & 53.43 & 732.37 & 0.24 & 729.61 & 7.10 & 1.36 & 46.08 & 1.09 & 0.31 & 0 & 0 & 402 & 287945 \\ 
& 5 & 53.01 & 1516.04 & 0.23 & 1512.93 & 7.10 & 1.37 & 45.68 & 1.09 & 0.65 & 0 & 0 & 402 & 359731 \\ \cline{ 2- 15}

\multirow{5}{*}{500} & 1 & 103.46 & 10.21 & 0.47 & 5.98 & 14.56 & 1.86 & 88.43 & 2.38 & 0 & 0 & 0 & 502 & 113418 \\ 
 & 2 & 103.23 & 87.44 & 0.42 & 83.17 & 14.56 & 1.85 & 88.25 & 2.37 & 0.05 & 0 & 0 & 502 & 225835 \\ 
 & 3 & 102.28 & 473.58 & 0.46 & 469.20 & 14.56 & 1.85 & 87.26 & 2.32 & 0.21 & 0 & 0 & 502 & 338252 \\ 
 & 4 & 102.60 & 1326.84 & 0.52 & 1322.02 & 14.56 & 1.85 & 87.51 & 2.34 & 0.63 & 0 & 0 & 502 & 450669 \\ 
 & 5 & 103.73 & 2589.89 & 0.50 & 2584.69 & 14.56 & 1.85 & 88.67 & 2.34 & 1.65 & 0.01 & 5 & 502 & 563086 \\ \cline{ 2- 15}
 
\multirow{5}{*}{700} & 1 & 318.71 & 56.62 & 1.60 & 50.02 & 59.45 & 3.23 & 257.66 & 3.36 & 0.01 & 0 & 0 & 702 & 219694 \\ 
& 2 & 318.78 & 391.64 & 1.52 & 384.99 & 59.45 & 3.23 & 257.81 & 3.34 & 0.07 & 0 & 0 & 702 & 437987 \\ 
& 3 & 321.05 & 1631.95 & 1.57 & 1624.82 & 59.45 & 3.24 & 260.03 & 3.34 & 0.81 & 0.01 & 1 & 702 & 656280 \\ 
& 4 & 320.33 & 3165.19 & 1.57 & 3158.33 & 59.45 & 3.23 & 259.31 & 3.38 & 1.56 & 0.12 & 9 & 702 & 874573 \\ 
& 5 & 319.01 & \multicolumn{1}{r}{\texttt{TL}} & 1.57 & \multicolumn{1}{r}{\texttt{TL}} & 59.45 & 3.25 & 258.00 & 3.34 & 0.03 & 28.21 & 10 & 702 & 1092866 \\ \cline{ 2- 15}

\multirow{5}{*}{920} & 1 & 937.45 & 548.75 & 14.23 & 535.54 & 216.81 & 5.05 & 706.40 & 8.01 & 0.15 & 0 & 0 & 922 & 370949 \\ 
& 2 & 927.92 & 2809.78 & 13.98 & 2795.74 & 216.81 & 5.06 & 697.13 & 7.99 & 1.03 & 0.12 & 2 & 922 & 740057 \\ 
& 3 & 927.02 & \multicolumn{1}{r}{\texttt{TL}} & 14.26 & \multicolumn{1}{r}{\texttt{TL}} & 216.81 & 5.06 & 695.94 & 8.16 & 0.18 & 0.94 & 10 & 922 & 1109165 \\ 
& 4 & 928.00 & \multicolumn{1}{r}{\texttt{TL}} & 14.01 & \multicolumn{1}{r}{\texttt{TL}} & 216.81 & 5.04 & 697.18 & 8.07 & 0.08 & 26.68 & 10 & 922 & 1478273 \\ 
& 5 & 932.15 & \multicolumn{1}{r}{\texttt{TL}} & 13.64 & \multicolumn{1}{r}{\texttt{TL}} & 216.81 & 5.05 & 701.70 & 7.93 & 0.01 & 38.09 & 10 & 922 & 1847381 \\ 
		\hline
		\end{tabular}}
		\caption{Performance of the computational experiments for $\rho(2)=0.005$}
		\label{t:computational1}
		\end{table}

\begin{table}[htbp]
	\centering
	\adjustbox{width=\textwidth}{\begin{tabular}{ll|rrrrrrrrr|rc|rr}
		\multicolumn{2}{l|}{} &\multicolumn{9}{c|}{\texttt{CPUTime (secs.)}} &\multicolumn{2}{c|}{} &\multicolumn{2}{c}{}\\
			
		\multicolumn{2}{l|}{} & \multicolumn{ 2}{c}{\texttt{Total}} & \multicolumn{ 2}{c}{\texttt{Solving}} & \multicolumn{ 2}{c}{\texttt{Prepr.}} & \multicolumn{ 2}{c}{\texttt{Ctrs. Gen.}} & \multicolumn{ 1}{c|}{\texttt{Callback}} & \multicolumn{ 1}{c}{\texttt{MIPGAP}} & \multicolumn{ 1}{c|}{\texttt{\#Unsolved}} & \multicolumn{ 2}{c}{\texttt{\#Ctrs}} \\  \cline{ 3- 15}
		
		n  & \multicolumn{1}{c|}{$p_2$} & \multicolumn{1}{c}{ \texttt{BIPS} } & \multicolumn{1}{c}{ \texttt{B\&C} } & \multicolumn{1}{c}{ \texttt{BIPS} } & \multicolumn{1}{c}{ \texttt{B\&C} } & \multicolumn{1}{c}{ \texttt{BIPS} } & \multicolumn{1}{c}{ \texttt{B\&C} } & \multicolumn{1}{c}{ \texttt{BIPS} } & \multicolumn{1}{c}{ \texttt{B\&C} } & \multicolumn{1}{c|}{ \texttt{B\&C} } & \multicolumn{1}{c}{ \texttt{B\&C} } & \multicolumn{1}{c|}{ \texttt{B\&C} } & \multicolumn{1}{c}{\texttt{BIPS}} & \multicolumn{1}{c}{ \texttt{B\&C} } \\ \hline	
\multirow{5}{*}{400} & 1 & 130.38 & 11.92 & 3.05 & 9.90 & 24.49 & 0.89 & 102.84 & 1.11 & 0.02 & 0 & 0 & 402 & 62344 \\ 
& 2 & 130.23 & 109.60 & 2.59 & 107.50 & 24.49 & 0.89 & 103.14 & 1.16 & 0.05 & 0 & 0 & 402 & 123887 \\ 
& 3 & 129.79 & 691.50 & 2.64 & 689.38 & 24.49 & 0.89 & 102.67 & 1.10 & 0.13 & 0 & 0 & 402 & 185430 \\ 
& 4 & 130.17 & 2252.92 & 2.87 & 2250.76 & 24.49 & 0.89 & 102.80 & 1.11 & 0.22 & 0.01 & 4 & 402 & 246973 \\ 
& 5 & 130.23 & 3299.38 & 2.98 & 3297.28 & 24.49 & 0.93 & 102.76 & 1.22 & 0.22 & 0.45 & 9 & 402 & 308516 \\ \cline{ 2- 15}
	
\multirow{5}{*}{500} & 1 & 288.18 & 22.46 & 9.13 & 18.71 & 66.61 & 1.35 & 212.44 & 2.36 & 0.04 & 0 & 0 & 502 & 95286 \\ 
 & 2 & 289.67 & 427.25 & 9.93 & 423.48 & 66.61 & 1.34 & 213.14 & 2.34 & 0.09 & 0 & 0 & 502 & 189571 \\ 
 & 3 & 289.07 & 1467.77 & 9.70 & 1463.84 & 66.61 & 1.34 & 212.76 & 2.34 & 0.26 & 0.01 & 1 & 502 & 283856 \\ 
 & 4 & 289.85 & 2614.43 & 9.89 & 2610.61 & 66.61 & 1.34 & 213.36 & 2.34 & 0.24 & 0.42 & 9 & 502 & 378141 \\ 
 & 5 & 289.29 & \multicolumn{1}{r}{\texttt{TL}} & 9.92 & \multicolumn{1}{r}{\texttt{TL}} & 66.61 & 1.34 & 212.76 & 2.36 & 0.03 & 28.37 & 10 & 502 & 472426 \\ \cline{ 2- 15}
 
\multirow{5}{*}{700} & 1 & 1001.96 & 67.20 & 40.67 & 61.42 & 325.38 & 2.35 & 635.91 & 3.36 & 0.08 & 0 & 0 & 702 & 180619 \\ 
& 2 & 1000.09 & 1519.78 & 41.36 & 1513.80 & 325.38 & 2.35 & 633.35 & 3.36 & 0.28 & 0 & 1 & 702 & 359837 \\ 
& 3 & 999.35 & 3289.82 & 40.67 & 3284.10 & 325.38 & 2.35 & 633.30 & 3.35 & 0.16 & 4.14 & 9 & 702 & 539055 \\ 
& 4 & 1000.46 & \multicolumn{1}{r}{\texttt{TL}} & 41.06 & \multicolumn{1}{r}{\texttt{TL}} & 325.38 & 2.35 & 634.02 & 3.33 & 0.01 & 30.58 & 10 & 702 & 718273 \\ 
& 5 & 999.80 & \multicolumn{1}{r}{\texttt{TL}} & 42.02 & \multicolumn{1}{r}{\texttt{TL}} & 325.38 & 2.35 & 632.40 & 3.35 & 0.01 & 36.89 & 10 & 702 & 897491 \\ \cline{ 2- 15}

\multirow{5}{*}{920} & 1 & 3253.39 & 471.14 & 141.97 & 458.74 & 1401.62 & 3.84 & 1709.81 & 7.94 & 0.63 & 0 & 0 & 922 & 292544 \\ 
& 2 & 3270.67 & 3320.87 & 141.63 & 3303.30 & 1401.62 & 3.84 & 1727.42 & 8.07 & 0.97 & 5.50 & 9 & 922 & 583247 \\ 
& 3 & 3267.82 & \multicolumn{1}{r}{\texttt{TL}} & 142.66 & \multicolumn{1}{r}{\texttt{TL}} & 1401.62 & 3.83 & 1723.55 & 7.95 & 0.01 & 22.84 & 10 & 922 & 873950 \\ 
& 4 & 3268.61 & \multicolumn{1}{r}{\texttt{TL}} & 145.35 & \multicolumn{1}{r}{\texttt{TL}} & 1401.62 & 3.85 & 1721.64 & 8.07 & 0.01 & 29.14 & 10 & 922 & 1164653 \\ 
& 5 & 3258.10 & \multicolumn{1}{r}{\texttt{TL}} & 147.56 & \multicolumn{1}{r}{\texttt{TL}} & 1401.62 & 3.83 & 1708.92 & 7.93 & 0.01 & 35.09 & 10 & 922 & 1455356 \\ \hline
		\end{tabular}}
		\caption{Performance of the computational experiments for $\rho(2)=0.01$}
		\label{t:computational2}
		\end{table}

\begin{figure}[h]
\begin{center}
\adjustbox{scale=0.82}{
	\begin{subfigure}[b]{.45\textwidth}
		\centering
		\includegraphics[scale=1.15]{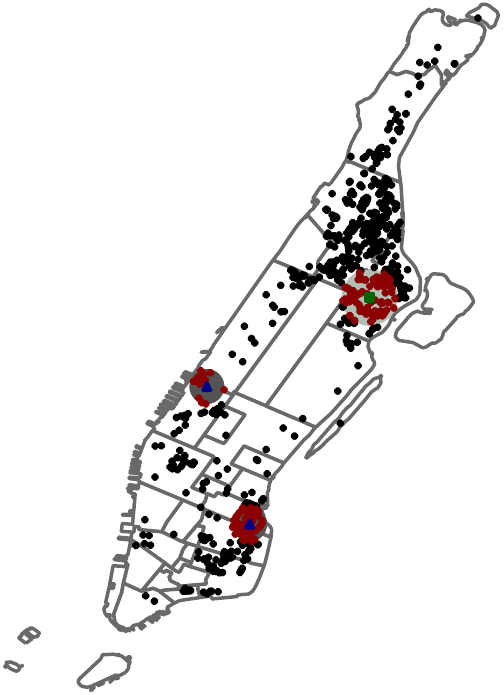}
		\caption{$p_1 = 1$ and $p_2 = 2$.}\label{f:mc1}
	\end{subfigure}~\begin{subfigure}[b]{.45\textwidth}
		\centering
		\includegraphics[scale=1.15]{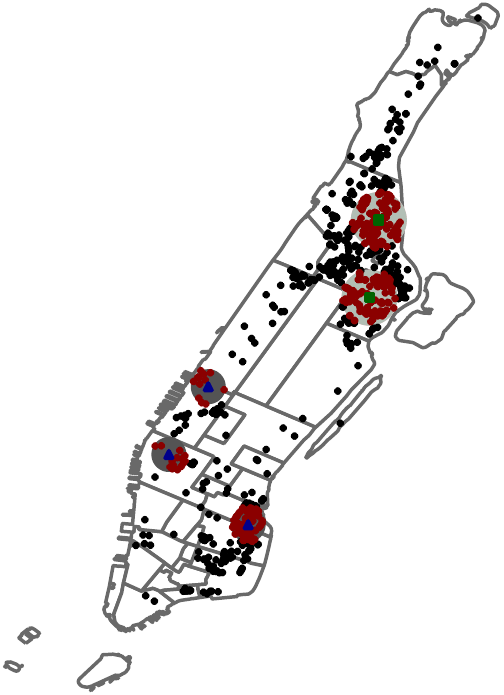}
		\caption{$p_1 = 2$ and $p_2 = 3$.}
	\end{subfigure}}\\ \quad \\
	\adjustbox{scale=0.82}{
	\begin{subfigure}[b]{.45\textwidth}
		\centering
		\includegraphics[scale=1.15]{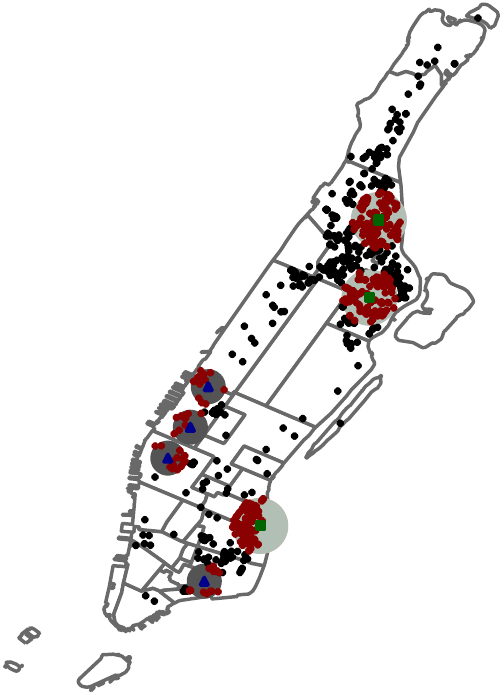}
		\caption{$p_1 = 3$ and $p_2 = 4$.}
	\end{subfigure}~
	\begin{subfigure}[b]{.45\textwidth}
		\centering
		\includegraphics[scale=1.15]{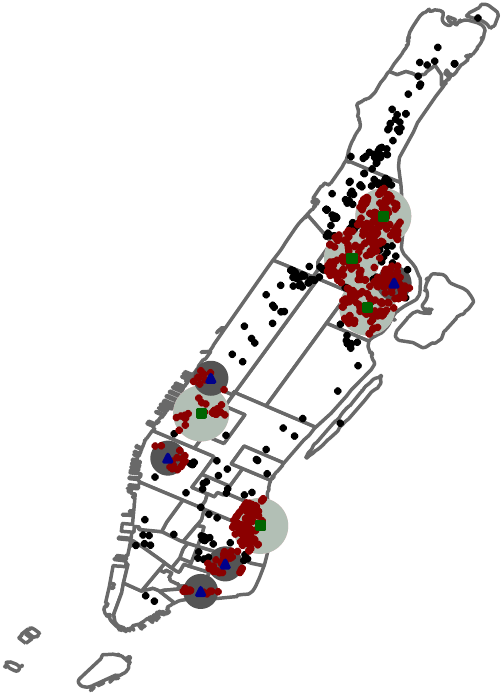}
		\caption{$p_1 = 5$ and $p_2 = 5$. }\label{f:mc4}
	\end{subfigure}}
	\caption{Solutions of some instances of our testbed for $n = 920$.}
\label{f:manhattan_cov}
\end{center}
\end{figure}

In Figure \ref{f:manhattan_cov}, we show the solutions of four instances of our testbed obtained with the \texttt{BIPS} approach painted on the Manhattan map. We draw the solutions for the larger instances ($n=920$ demand points), radii $\rho(1) = 0.008$ and $ \rho(2) = 0.005$ and different values of $p_1$ and $p_2$. In the figure, red dots represent the covered demand nodes, green squares the positions of the discrete facilities, and blue triangles the positions of continuous facilities; the coverage areas for the discrete facilities are drawn in light green color and those of the continuous facilities are colored in gray. The percentages of covered demand for these instances ranges between $29\%$ (Figure \ref{f:mc1}) and $75\%$ (Figure \ref{f:mc4}).

\section{Conclusions and Further Research}

In this paper we investigated the multi-type maximal covering facility location problem.
A general modeling framework was discussed, which was adapted to an hybridized discrete-continuous facility location problem.
In this case we could go deeper in our analysis.
For the particular case in which the space underlying the continuous location problem is the Euclidean space and when the euclidean norm is used a third model could be proposed.
The results highlighted the viability of a branch-and-cut algorithm for dealing with the problem in its general form. 
In particular, instances with up to 920 demand nodes and two types of facilities (discrete and continuous) could be solved rather efficiently. 
This defines a new state-of-the-art in terms of maximal covering location problems with a large potential number of locations for the discrete facilities.

The work done encourages some other research lines.
These include more work on the development of valid inequalities for the general integer linear programming model thus leading to an even better polyhedral description of the feasibility set, the inclusion of time as a decision dimension, and the inclusion of uncertainty either in the demand or in the number of facilities that can be open.

\section*{Acknowledgements}

This research has been partially supported by Spanish Ministerio de Ciencia e Innovación, AEI/FEDER grant number PID2020-114594GBC21, Junta de Andalucía projects P18-FR-1422/2369 and projects FEDERUS-1256951, B-FQM-322-UGR20, CEI-3-FQM331 and NetmeetData (Fundación BBVA 2019). The first author was also partially supported by the IMAG-Maria de Maeztu grant CEX2020-001105-M /AEI /10.13039/501100011033.
The second author was partially supported by Spanish Ministry of Education and Science grant number PEJ2018-002962-A, the PhD Program in Mathematics at the Universidad de Granada and Becas de Movilidad entre Universidades Andaluzas e Iberoamericanas (AUIP).
The third author was partially funded by grant UIDB/04561/2020 from National Funding from FCT---Fundação para a Ciência e Tecnologia, Portugal.

%\bibliographystyle{elsarticle-harv} 
%\bibliography{references}

\end{document}